\documentclass{amsart}
\usepackage{amsmath}
\usepackage{graphicx}
\usepackage{latexsym}
\usepackage{amsfonts}
\usepackage{amssymb}
\usepackage{graphicx}
\usepackage{dsfont}
\usepackage[utf8]{inputenc}
\usepackage{wasysym}
\usepackage{xfrac}
\usepackage{float}
\usepackage{hyperref}
\usepackage[hyphenbreaks]{breakurl}
\usepackage{tikz}

\usepackage{enumitem}
\setlist{nosep}

\usepackage[shortcuts]{extdash}

\usepackage{romannum}

\usepackage{todonotes}

\usepackage{mathtools}
\usepackage{bbm}
\usepackage{accents}

\allowdisplaybreaks
\usepackage[style=numeric,backend=biber,maxnames=99,isbn=false]{biblatex} 
\defbibheading{secbib}[\bibname]{%
  \section*{#1}%
  \markboth{#1}{#1}}
  
\renewbibmacro*{url}{%
  \iftoggle{bbx:url}
    {\iffieldundef{doi}{\printfield{url}}{}}
    {}%
  \newunit\newblock
}

\theoremstyle{plain}
\newtheorem{theorem}{Theorem}[section]
\newtheorem*{theorem*}{Theorem}

\newtheorem{lemma}[theorem]{Lemma}
\newtheorem{proposition}[theorem]{Proposition}
\theoremstyle{definition}

\newtheorem{definition}[theorem]{Definition}
\newtheorem{example}[theorem]{Example}

\newtheorem{remark}[theorem]{Remark}
\newtheorem*{remark*}{Remark}

\newcommand{\N}{\mathbb{N}}

\newcommand{\R}{\mathbb{R}}

\newcommand{\Abs}[1]{\left\vert #1 \right\vert}

\newcommand{\abs}[1]{\vert #1 \vert}
\newcommand{\If}{\mathbbm{1}}
\newcommand{\IF}[1]{\mathbbm{1}_{\{ #1 \}}}
\newcommand{\sgn}{\operatorname{sgn}}
\newcommand{\supp}{\operatorname{supp}}
\newcommand{\suppex}{\overline{\operatorname{supp}}}

\newcommand{\EW}[1]{\mathbb{E}\left[ #1 \right]}

\newcommand{\EWx}[2]{\mathbb{E}_{#1}\left[ #2 \right]}

\renewcommand{\Pr}{\mathbb{P}}
\newcommand{\PR}[1]{\mathbb{P}\left( #1 \right)}

\newcommand{\PRx}[2]{\mathbb{P}_{#1} \left( #2 \right)}

\newcommand{\PRb}[2]{\mathbb{P}\left( #1 \; \middle| #2 \right)}
\newcommand{\PRbx}[3]{\mathbb{P}_{#1} \left( #2 \; \middle| #3 \right)}

\newcommand{\st}{\preceq_{\operatorname{st}}}
\newcommand{\hr}{\preceq_{\operatorname{hr}}}

\newcommand{\Hconv}{\overset{\Gamma}{\to}}

\newcommand{\intdiff}{\mathop{}\!\mathrm{d}}

\newcommand{\epi}{\operatorname{epi}}

\begin{document}

\pagenumbering{arabic}
\title[Existence and uniqueness of the inverse first-passage time problem]{Conditions for existence and uniqueness\\ of the inverse first-passage time problem\\ applicable for Lévy processes and diffusions}

\author[Klump]{Alexander Klump$^*$}
\address[*,\dag]{Institute of Mathematics and Informatics, Bulgarian Academy of Sciences, Akad. Georgi Bonchev str., bl. 8, Sofia 1113, Bulgaria}
\email{aklump@math.upb.de}
\author[Savov]{Mladen Savov$^\dag$}
\address[*]{Faculty of Mathematics and Informatics, Sofia University ”St. Kliment Ohridski”, bl. James Bouchier 5, Sofia 1164, Bulgaria}
\email{msavov@fmi.uni-sofia.bg}

\subjclass{60G07, 60J25, 60G51, 60J60; 60G40, 60E15}
\keywords{inverse first-passage time problem, stochastic process, Markov process, Lévy process, diffusion}
        
\begin{abstract}
For a stochastic process $(X_t)_{t\geq0}$ we establish conditions under which the inverse first-passage time problem has a solution for any random variable $\xi >0$. For Markov processes we give additional conditions under which the solutions are unique and solutions corresponding to ordered initial states fulfill a comparison principle. As examples we show that these conditions include Lévy processes with infinite activity or unbounded variation and diffusions on an interval with appropriate behavior at the boundaries. Our methods are based on the techniques used in the case of Brownian motion and rely on discrete approximations of solutions via $\Gamma$-convergence from \cite{Anu80} and \cite{chen2011} combined with stochastic ordering arguments adapted from \cite{kk2022}.
\end{abstract}
\maketitle

\section{Introduction}
Given a random variable $\xi$ with values in $(0,\infty )$ the inverse first-passage time problem for a stochastic process $(X_t)_{t\geq0}$ with values in $\R$ consists of the question whether there exists $b:[0,\infty] \to [-\infty,\infty]$ such that the first-passage time
\begin{align*}
\tau_b \coloneqq \inf\{t> 0 : X_t \geq b(t)\}
\end{align*}
of $b$ has the same distribution as $\xi$. If this question can be answered affirmatively, one naturally asks whether these solutions are unique in a reasonable sense and which properties they have. The terminology is due to the first-passage time problem, where for a stochastic process $(X_t)_{t\geq0}$ and a function $b$ the question is to determine properties of the distribution of $\tau_b$. Primarily, the inverse first-passage time problem was studied for Brownian motion and revealed relations to free boundary problems, integral equations and optimal stopping problems and gave rise to applications in mathematical finance. For general processes this problem is of particular interest due to its possible relevance in applications and the new theoretical questions it gives rise to.

The inverse first-passage time problem roots back to a question of Shiryaev, whether there is a stopping time with respect to a Brownian motion which is exponentially distributed. This question was answered by \cite{dudley1977} for a general stochastic process by establishing conditions under which stopping times  with given distributions exist. For the inverse first-passage time problem the existence of lower-semicontinuous solutions was established in \cite{Anu80} for reflected Brownian motion by a discrete approximation of their epigraphs. In the case that $\xi >0$ has no atoms, existence and uniqueness have been obtained for diffusions in \cite{cheng2006} and \cite{chen2011} via a transfer into a free boundary problem. For Brownian motion uniqueness was shown for arbitrary $\xi >0$ in \cite{Eks16} via discretization of a related optimal stopping problem and independently deduced in \cite{beiglboeck2018} in a more general setting of optimal stopping problems with distribution constraints. For reflected Brownian motion, in \cite{kk2022} the uniqueness was shown via a discrete approximation argument paired with stochastic ordering. The discretizations of \cite{Eks16} and \cite{kk2022} are related to the approximation in \cite{Anu80}. Conditions for continuity of solutions were given in \cite{chen2011}, \cite{Eks16} and \cite{potiron2021}, where higher order regularity was studied in \cite{chen2022}. The inverse first-passage time problem and the first-passage time problem for Brownian motion are related to certain integral equations, see \cite{peskir2002b}, \cite{cheng2006}, \cite{Eks16}, \cite{kreinin2009a}. Numerical approaches for the case of Brownian motion are to be found in \cite{zucca2009}, \cite{abundo2006}, \cite{song2011}, \cite{poetzelberger2021}, \cite{k2023} and for an Ornstein-Uhlenbeck process in \cite{civallero2019}. Applications have been proposed in \cite{hull2001}, \cite{avellaneda2001}, \cite{sacerdote2005}, \cite{sacerdote2007}, \cite{yingjun2014} in the context of modeling default, neuronal activity or failure. For Brownian motion the modification of the problem to fix both $b$ and $\xi$ and to ask whether $X_0$ can be randomly distributed such that $\tau_b$ has distribution according to $\xi$ has been studied in \cite{kreinin2009b}, \cite{kreinin2009c}, \cite{abundo2013}, \cite{kreinin2014} and is naturally related to our comparison principle. Moreover, the inverse first-passage time problem for the case of Brownian motion and exponentially  distributed $\xi$ is related to \cite{demasi2019}, \cite{berestycki2019}, \cite{lee2020}, \cite{berestycki2022} and \cite{berestycki2021}, where hydrodynamic limits of  certain particle systems and corresponding free boundary problems are studied. For general $\xi$ a related particle system whose hydrodynamic limit is characterized by the inverse first-passage time problem for reflected Brownian motion has been constructed in \cite{k2023}. Another so-called soft-killing variant of the problem asks the same question but with a smoothed-out version of $\tau_b$, where one additionally waits for an exponential clock to ring after passing the boundary, and was treated in \cite{hening2014}, \cite{hening2020} and \cite{kk2022sk}. For a more detailed overview of related work in the case of Brownian motion we refer to \cite{k2022}. Another modification of the problem is to fix $b$ and the distributions of $\xi$ and $X_0$ and to search for a suitable stochastic process $(X_t)_{t\geq 0}$ in order to achieve that $\tau_b$ has the same distribution as $\xi$. This has been studied in \cite{davis2010} for a certain family of It\^o diffusions and in \cite{davis2015} for a family of processes obtained from deterministically time-changing a fixed Lévy process.

Let us summarize the methods used in this present paper. For the existence of solutions we pick up the idea of \cite{Anu80} from the case of reflected Brownian motion. A main ingredient in the proof of \cite{Anu80} is the continuity of the paths. By a careful adaption of the proof and by utilizing facts about left and right discontinuity of arbitrary functions we are able to work with quasi-left-continuous cádlág paths instead. For the uniqueness of solutions the Markov property allows to work with stochastic orders as in \cite{kk2022}. In this situation we present a new and elementary argument how to infer uniqueness from the discretization of \cite{Anu80}. Instead of using the Wasserstein distance for the marginal distribution as in \cite{kk2022} we use an adapted approximation of lower semicontinuous functions from \cite{chen2011} and \cite{Eks16}.

We want to emphasize the following relations to our results.
\begin{itemize}
\item For Brownian motion relations to free boundary problems \cite{chen2011}, \cite{berestycki2019}, optimal stopping problems \cite{Eks16}, \cite{beiglboeck2018}, integral equations \cite{kreinin2009a}, \cite{peskir2002b} and particle systems \cite{demasi2019}, \cite{k2023} are known. The question arises whether such relations extend to other stochastic processes.
\item By the comparison principle, the inverse first-passage time problem is related to the first-passage time problem and the modified problems studied in \cite{kreinin2009b} or \cite{davis2015}, which therefore gain interest in this general setting.
\item Continuity of the paths was also a main ingredient in the soft-killing inverse first-passage time problem for Brownian motion in \cite{hening2020}, \cite{kk2022sk}. It is natural to ask if this can also be extended to processes with discontinuous paths.
\item Our results give more credibility to numerical approaches where existence and uniqueness were assumed for the underlying processes, see \cite{zucca2009}, \cite{civallero2019}.
\item The existence result allows more flexibility in applications, since the process used for modeling purposes as in \cite{hull2001}, \cite{sacerdote2005}, \cite{yingjun2014} can be chosen from a broader range.
\item The discrete approximation of our approach yields a possible numerical Monte-Carlo type approximation similar as in \cite{k2023} for Brownian motion.
\end{itemize}

The paper is organized as follows. In Section~\ref{s:notation main result} we present our main results regarding existence, uniqueness and comparison principle as well as the conditions for Lévy processes and diffusions. The proofs of the main results are to be found in Section~\ref{s:existence}, Section~\ref{s:uniqueness} and Section~\ref{section:comparison principle}. The proofs regarding the conditions for Lévy processes as well as for diffusions are contained in Section~\ref{conditions for Levy} and Section~\ref{conditions for diffusions}.

\section{Main results}\label{s:notation main result}

\begin{definition}
We call a lower semicontinuous function $b:[0,\infty] \to [-\infty , \infty]$ \emph{boundary function}.
\end{definition}

Let $(\Omega , \mathcal{F}, \Pr )$ be a probability space endowed with a filtration $(\mathcal{F}_t)_{t\geq 0}$ fulfilling the usual conditions. Let $(X_t)_{t\geq 0}$ be an adapted stochastic process. For a boundary function $b$ define, in addition to $\tau_b$, the first-passage time variant
\begin{align*}
\tau_{b}' &\coloneqq \inf\{ t > 0 : X_t > b(t) \}.
\end{align*}

\begin{definition}
We say that $(X_t)_{t\geq 0}$ is \emph{quasi-left-continuous}, if for any non-decreasing sequence of stopping times $T_1\leq T_2 \leq \ldots$ and $T\coloneqq \lim_{n\to\infty} T_n$ it holds that
\begin{align*}
\lim_{n\to\infty} X_{T_n} = X_T
\end{align*}
almost surely on $ \{T<\infty\}$.
\end{definition}
\noindent For the existence of solutions the following assumptions are put in force:
\begin{enumerate}[label=(E\arabic*)]
\item\label{existence: diffusitivy assumption} For every $t>0$ the probability measure $ \PR{X_t \in \cdot\;}$ is diffuse.
\item\label{existence: continuity assumption} $(X_t)_{t\geq 0}$ has $\Pr$-a.s. right-continuous paths and is quasi-left-continuous.
\item\label{existence: first-passage time assumption} For any boundary function $b$ it holds $\Pr$-a.s. that $\tau_b = \tau_b'$.
\end{enumerate}

\begin{theorem}[Existence]\label{existence}
Assume that \ref{existence: diffusitivy assumption}, \ref{existence: continuity assumption} and \ref{existence: first-passage time assumption} are fulfilled.\\ Then, given a random variable $\xi >0$, there exists a boundary function $b$ such that
\begin{equation*}
\tau_b \overset{\text{d}}{=} \xi.
\end{equation*}
\end{theorem}

\begin{remark}
Let us emphasize that the conditions for existence are not very restrictive. Condition \ref{existence: diffusitivy assumption} is necessary to have solutions for any $\xi >0$. Furthermore, although conditions \ref{existence: continuity assumption} and \ref{existence: first-passage time assumption} are primarily employed for technical reasons, they turn out to be very natural. For example, \ref{existence: continuity assumption} is fulfilled by Feller processes with cádlág paths. A condition as \ref{existence: first-passage time assumption} is needed to exclude behavior as in Example~\ref{example}.
\end{remark}

\begin{remark}
If the process takes values in an interval $(0,R]$ with $R< \infty$, then the assumption of \ref{existence: first-passage time assumption} implies, by taking $b\equiv R$, that $\PR{\tau_R < \infty} = 0$.
\end{remark}

\begin{example}\label{example}
Let $X_t = - \abs{B_t + 1}$, where $(B_t)_{t\geq 0}$ is a standard Brownian motion. Then $\tau_0$ has the Lévy distribution with scale $1$, in particular it is supported on $(0,\infty )$. On the other hand, for every $\xi >0 $ which is supported on $(0,\infty )$ a solution $b$ must have values in $(-\infty , 0]$. This means that such a $\xi >0$ whose distribution is strictly larger in the usual stochastic order than the distribution of $\tau_0$ (Definition~\ref{usual stochastic order}) cannot be realized as first-passage time.
\end{example}

\begin{definition}\label{usual stochastic order}
For two probability measures $\mu , \nu$ on $(\R , \mathcal{B}(\R))$, we say $\mu$ is \emph{smaller in the usual stochastic order}, write $\mu \st \nu$, if and only if
\begin{align*}
\mu ((-\infty , c] )  \geq \nu ((-\infty , c]) \qquad \forall c \in \R.
\end{align*}
\end{definition}
\noindent For a measure $\mu$ on $\R$ we define its support as
\begin{align*}
\supp (\mu) \coloneqq \{x\in \R : \mu (U) > 0 \text{ whenever } U\subseteq \R \text{ is open and } x\in U \}.
\end{align*}

For a random variable $\xi >0$ we define
\begin{align*}
t^\xi \coloneqq \sup\supp\left( \PR{ \xi \in \cdot\;}\right) = \sup\{t > 0 : \PR{\xi > t} > 0\}.
\end{align*}
For a process taking values in an interval $E\subseteq \R$ the following assumptions are put in force for uniqueness:
\begin{enumerate}[label=(U\arabic*)]
\item\label{uniqueness: markov assumption} There is a family of probability measures $(\Pr_x)_{x\in E}$ such that
\begin{align*}
\PRx{x}{X_t \in E \;\forall t \geq 0} = 1 \quad\forall x\in E
\end{align*}
and $((X_t)_{t\geq 0}, ( \mathcal{F}_t)_{t\geq 0} , (\Pr_x)_{x\in E})$ is a Markov process as per \cite[Vol.1,p.77]{dynkin1965}.
\item\label{uniqueness: order preservation assumption} For probability measures $\mu_1,\mu_2$ on $E$ and all $t>0$, we have that $\mu_1 \st \mu_2$ implies
\begin{align*}
\PRx{\mu_1 }{X_t \in \cdot \;} \st \PRx{\mu_2}{X_t \in \cdot\;}.
\end{align*}
\item\label{uniqueness: support assumption} For a random variable $\xi >0$ there is $I^\xi \subseteq (0,t^\xi )$ such that for any boundary function $b$ with values in $\overline{E}$ and $\tau_b \overset{\text{d}}{=} \xi$ we have
\begin{align*}
b(t) = \sup \supp \left(\PRb{X_t \in \cdot}{\tau_{b} > t}\right)
\end{align*}
for all $t \in I^\xi$.
\end{enumerate}
For a probability measure $\mu$ on $E$ we define $\Pr_\mu \coloneqq \int_E \Pr_x \mu (\intdiff x)$.

\begin{theorem}[Uniqueness]\label{uniqueness}
Let $E \subseteq \R$ be an interval and denote $\overline{E} = [L , R ]$ with $L,R\in [-\infty, \infty]$. Fix a probability measure $\mu$ on $E$. Let $\xi >0$ be a random variable. Assume that \ref{uniqueness: markov assumption}, \ref{uniqueness: order preservation assumption} are fulfilled and \ref{existence: diffusitivy assumption}, \ref{existence: continuity assumption}, \ref{existence: first-passage time assumption} and \ref{uniqueness: support assumption} with $\Pr \coloneqq \Pr_\mu$ and $I^\xi\subseteq (0,t^\xi )$ are fulfilled.\\ Then all boundary functions $b$ with values in $\overline{E}$ and $\tau_b \overset{\text{d}}{=} \xi$ under $\Pr$ coincide on $ I^\xi$.
\end{theorem}

\begin{remark}
Let us comment on the conditions for uniqueness. \ref{uniqueness: markov assumption} and \ref{uniqueness: order preservation assumption} are contingent on our method of proof, for which we do not anticipate problems when working with inhomogeneous Markov processes instead. However, for simplicity we refrain from doing this. \ref{uniqueness: support assumption} is necessary for uniqueness, since otherwise we could alter values of a solution $b$ on the set $I^\xi$ without affecting the distribution of the first-passage time $\tau_b$.
\end{remark}

\begin{definition}
For two random variables $\xi ,\zeta >0$ we say \emph{$\xi$ is smaller in the hazard rate order than $\zeta$}, write $\xi \hr \zeta$, if
\begin{align*}
[0,t^\zeta) \to [0,1] ,\; t \mapsto \frac{\PR{\xi > t}}{ \PR{\zeta > t}}
\end{align*}
is a non-increasing function.
\end{definition}

\begin{theorem}[Comparison principle]\label{comparison principle}
Let $E \subseteq \R$ be an interval. Fix two probability measures $\mu_1 , \mu_2$ on $E$ such that $\mu_1 \st \mu_2$. Let $\xi_1, \xi_2 >0$ be random variables such that $\xi_1 \hr \xi_2$.  Assume that \ref{uniqueness: markov assumption}, \ref{uniqueness: order preservation assumption} are fulfilled and \ref{existence: diffusitivy assumption}, \ref{existence: continuity assumption}, \ref{existence: first-passage time assumption} with $\Pr \coloneqq \Pr_{\mu_i}$ are fulfilled.\\ Then for $i\in\{1,2\}$ there exist boundary functions $b_i$ with $\tau_{b_i} \overset{\text{d}}{=} \xi_i$ under $\Pr_{\mu_i}$ such that
\begin{align*}
b_{1} \leq b_{2}
\end{align*}
pointwise.
\end{theorem}

\begin{remark}
Since in Theorem~\ref{comparison principle} we merely state the existence of ordered solutions, we can spare \ref{uniqueness: support assumption} in our list of assumptions.
\end{remark}

\subsection{Lévy processes}

In Section~\ref{conditions for Levy} we establish conditions for Lévy processes under which we can apply Theorem~\ref{existence} and Theorem~\ref{uniqueness}. We will summarize these conditions below in Theorem~\ref{Levy processes result}.

We say a Lévy process has characteristic a triple $(a,\sigma^2 , \Pi )$ if
\begin{align}\label{characteristic triple}
- \log \left( \EW{\exp (i\theta X_1) } \right) = i\theta a + \frac{\sigma^2}{2} \theta^2 + \int_\R ( 1 - e^{i\theta x} + i\theta x 1_{(-1,1)}(x))\Pi (\intdiff x),
\end{align}
where $a\in \R$, $\sigma^2 \geq 0$ and $\Pi$ is a measure on $\R \setminus \{0\}$ with $\int_\R (1\wedge x^2) \Pi (\intdiff x) < \infty$. If $(X_t)_{t\geq 0}$ is a Lévy process with $\PR{X_0 = 0} = 1$, then for $x\in \R$ let $\Pr_x$ be a measure such that $\PRx x{(X_t )_{t\geq 0} \in \cdot\;}\coloneqq \PR{ (X_t +x)_{t\geq 0} \in \cdot\;}$. For the following statement note that $\PRx{0}{X_1 \in \cdot\;}$ is diffuse if and only if we have that $\sigma^2 >0$ or $\Pi (\R) = \infty$, see for instance Theorem~27.4 in \cite{Sato1999}. Equivalently, we could say the Lévy process is not a compound Poisson process with or without drift. On the other hand, if a probability measure $\mu$ on $\R$ is diffuse, then $\PRx \mu {X_t\in\cdot\;}$ is diffuse.

\begin{theorem}[Lévy processes]\label{Levy processes result}
Let $(X_t)_{t\geq 0}$ be a Lévy process with characteristic triple $(a,\sigma^2 , \Pi)$, $\xi>0$ be a random variable and $\mu$ be a probability measure on $\R$. Then \ref{uniqueness: markov assumption} and \ref{uniqueness: order preservation assumption} are fulfilled with $E=\R$ and \ref{existence: continuity assumption} is fulfilled with $\Pr \coloneqq \Pr_{\mu}$.
\paragraph{Existence:}
We have that \ref{existence: diffusitivy assumption} implies \ref{existence: first-passage time assumption}. In particular, assuming \ref{existence: diffusitivy assumption}, there exists a boundary function $b$ such that $\tau_b \overset{\text d} = \xi$ under $\Pr$.
\paragraph{Uniqueness:}
\begin{enumerate}
\item[(a)] Let one of the following be fulfilled:
\begin{itemize}
\item[(a.i)] $(X_t)_{t\geq 0}$ has unbounded variation, i.e. $\sigma^2 > 0$ or $\int_\R (1\wedge \abs{x}) \Pi (\intdiff x) = \infty$,
\item[(a.ii)] $0\in\supp (\Pi)$ and $\Pi ((0,\infty )) > 0$.
\end{itemize}
Then \ref{uniqueness: support assumption} is fulfilled with $I^\xi \coloneqq (0,t^\xi )$.
\item[(b)] Let $0\in\supp (\Pi)$ and $\Pi ((-\infty , 0 )) > 0$.\\ Then \ref{uniqueness: support assumption} is fulfilled with $I^\xi \coloneqq \supp ( \PR{\xi \in \cdot\; })\cap (0,t^\xi )$.
\end{enumerate}
In particular, assuming \ref{existence: diffusitivy assumption} and ((a) or (b)), all boundary functions $b$ with $\tau_b \overset{\text d} = \xi$ under $\Pr$ coincide on $I^\xi$.
\end{theorem}

\begin{remark}
In order to demonstrate the phrasing of Theorem~\ref{Levy processes result} let us mention some examples of Lévy processes:
\begin{itemize}
\item If the Lévy measure has infinite activity, i.e. $\Pi (\R) = \infty$, then we have $0 \in \supp (\Pi)$ and \ref{existence: diffusitivy assumption}, and thus we have existence and uniqueness.
\item If $(X_t)_{t\geq 0}$ has a Brownian component we have existence and uniqueness.
\item If the law of $X_0$ is diffuse and $(X_t)_{t\geq 0}$ is a Poisson process with jumps of constant height we have existence.
\item If $(-X_t)_{t\geq 0}$ is a Gamma process and $\xi \sim \text{Exp}$, we have existence and uniqueness of solutions on $(0,\infty)$.
\end{itemize}
\end{remark}

\subsection{Diffusions on an interval}

In Section~\ref{conditions for diffusions} we establish conditions for diffusions under which we can apply Theorem~\ref{existence} and Theorem~\ref{uniqueness}. We will summarize these conditions below in Theorem~\ref{diffusions result}.

For the definition of a diffusion on an interval we adapt Definition~5.20 of \cite{karatzas1991}.
\begin{definition}\label{definition diffusion}
Let $E\subseteq \R$ be an interval and $\overline{E} = [L,R] \subseteq [-\infty , \infty]$. Furthermore, let $\beta : E \to \R $ and $\sigma : E \to \R$ be Borel-measurable functions. Let $(\Pr_x)_{x\in E}$ be a family of probability measures and $(X_t)_{t\geq 0}$, $(B_t)_{t\geq 0}$ stochastic processes such that
\begin{itemize}
\item[(i)] $((X_t)_{t\geq 0} , (\mathcal{F}_t)_{t\geq 0}, (\Pr_x)_{x\in E})$ is a strong Markov process on $E$,
\item[(ii)] $(X_t)_{t\geq 0}$ is $(\mathcal{F}_t)_{t\geq 0}$-adapted with continuous paths $\Pr_x$-a.s. and $(B_t)_{t\geq 0}$ is a Brownian motion with respect to $(\mathcal{F}_t)_{t\geq 0}$ and $\Pr_x$ for every $x\in E$,
\item[(iii)] with strictly monotone sequences $(\ell_n)_{n\in \N}$ and $(r_n)_{n\in\N}$  satisfying $L < \ell_n < r_n < R$, $\lim_{n\to\infty } \ell_n = L$ and $\lim_{n\to\infty} r_n = R$ and
\begin{align*}
S_n \coloneqq \inf \{t\geq 0 : X_t \notin (\ell_n , r_n) \}
\end{align*}
it holds for all $t\geq 0$ that
\begin{align*}
\PRx{x}{ \int_0^{t\wedge S_n} \abs{\beta (X_s)} + \sigma^2 (X_s) \intdiff s <\infty} = 1
\end{align*}
and
\begin{align}\label{diffusion equation}
\PRx{x}{X_{t\wedge S_n } = x+ \int_0^{t\wedge S_n} \beta(X_s ) \intdiff s + \int_0^{t\wedge S_n} \sigma (X_s )\intdiff B_s \; \forall t \geq 0} =1
\end{align}
for all $n\in\N$ and all $x\in (L,R)$.
\end{itemize}
We call $(X_t)_{t\geq 0}$ a diffusion on $E$ with coefficients $\beta$ and $\sigma$.
\end{definition}

For the following statement, note that if $\PR{X_t\in\cdot \;}$ is diffuse, the process can hit the lower boundary $L$ in finite time, but it cannot get stuck there.

\begin{theorem}\label{diffusions result}
Let $(X_t)_{t\geq 0}$ be a diffusion on $E$ with coefficients $\beta$ and $\sigma$ such that $\sigma \in C^1 ((L,R))$, $\sigma >0$ and $\beta$ is locally bounded on $(L,R)$. Let $\xi >0$ be a random variable and $\mu$ be a probability measure on $E$. If $R\notin E$, then \ref{uniqueness: markov assumption}, \ref{uniqueness: order preservation assumption} are fulfilled and \ref{existence: continuity assumption}, \ref{existence: first-passage time assumption} are fulfilled with $\Pr \coloneqq \Pr_\mu$. If additionally \ref{existence: diffusitivy assumption} holds then we have \ref{uniqueness: support assumption} with $I^\xi = (0,t^\xi)$.
\paragraph{Existence and uniqueness:} In particular, assuming $R\notin E$ and \ref{existence: diffusitivy assumption}, there exists a unique boundary function $b$ on $(0,t^\xi)$ with values in $\overline E$ such that $\tau_b \overset{\text d} = \xi$.
\end{theorem}

\begin{remark}
In order to demonstrate the phrasing of Theorem~\ref{diffusions result} and the usage of the notion of diffusion on an interval, let us mention the following example. If $(X_t)_{t\geq 0}$ is a reflected Bessel process of dimension $\delta > 0$ on $E=[0,\infty )$ (cf. \cite[p.29]{lawler2019}), for $t < S = \inf\{ s \geq  0 : X_s = 0\}$ the process fulfills the stochastic integral equation
\begin{align*}
X_t = X_0 + \int_0^t \frac{\delta - 1}{2 X_s} \intdiff s +  B_t.
\end{align*}
The coefficients fulfill the conditions of Theorem~\ref{diffusions result} on $(0,\infty)$. Note that we have $R=\infty$. Since $\delta > 0$ the law of $X_t$ is absolutely continuous w.r.t. Lebesgue measure and thus \ref{existence: diffusitivy assumption} is fulfilled. This means that we have existence and uniqueness of solutions for any $\xi >0$.
\end{remark}

\section{Existence: Proof of Theorem~\ref{existence}}\label{s:existence}

Let us explain the role of the conditions in the proof of Theorem~\ref{existence}, which will give us a common thread. Condition \ref{existence: diffusitivy assumption} is necessary to have solutions for any $\xi >0$ and allows for the construction of a discrete approximation. For this approximation we will use the following notion of $\Gamma$-convergence. The conditions \ref{existence: continuity assumption} and \ref{existence: first-passage time assumption} ensure that this approximation provides a solution. The idea of this proof follows the approach of \cite{Anu80}.

\begin{definition}
We call a sequence $(b_n)_{n\in\N}$ of boundary functions $\Gamma$-convergent to a boundary function $b$, write $b_n \Hconv b$, if and only if
\begin{enumerate}
\item[(i)] for every convergent sequence $(t_n)_{n\in\N} \subset [0,\infty]$ with $\lim_{n\to\infty}t_n =t$ holds
\begin{align*}
\liminf_{n\to\infty} b_n (t_n) \geq b(t),
\end{align*}
\item[(ii)] for any $t\in [0,\infty]$ exists a convergent sequence $(t_n)_{n\in\N} \subset [0,\infty]$ with $\lim_{n\to\infty}t_n = t$ such that
\begin{align*}
\lim_{n\to\infty} b_n (t_n ) = b(t).
\end{align*}
\end{enumerate}
\end{definition}

For a boundary function $b$ and $s\geq 0$ we define
\begin{align*}
b|_s (t) \coloneqq \infty  \If_{[0,s)} (t) + b \If_{[s,\infty]}(t).
\end{align*}

\begin{proposition}\label{Anulova convergence}
Assume that $(X_t)_{t\geq 0}$ has a.s. right-continuous paths and is quasi-left-continuous, i.e. \ref{existence: continuity assumption} holds. Let $b$ be a boundary function and assume $\tau_{b\vert_s} \overset{\text{d}}{=} \tau_{b\vert_s}'$ for every $s>0$. Further let $b_n \Hconv b$ and assume that
\begin{align*}
\lim_{s\searrow 0}\limsup_{n\to\infty} \PR{\tau_{b_n} \leq s} =0.
\end{align*}
Then
\begin{align*}
\tau_{b_n} \overset{\Pr}{\to} \tau_b
\end{align*}
in probability as $n\to\infty$.
\end{proposition}

Proposition~\ref{Anulova convergence} will be proved after a sequence of preliminary lemmas.

\begin{remark}\label{Anulova convergence remark}
Note that, if $\xi >0$ is a random variable and $\tau_{b_n} \to \xi$ in distribution, then it follows by Portmanteau's theorem that
\begin{align*}
\lim_{s\searrow 0} \limsup_{n\to\infty} \PR{\tau_{b_n} \leq s} \leq \lim_{s\searrow 0} \PR{\xi \leq s} = 0.
\end{align*}
\end{remark}

\begin{remark}
For a boundary function $b$ we can rewrite $\tau_b$ as
\begin{align*}
\tau_b = \inf\{t >0 : (t,X_t) \in \epi (b)\},
\end{align*}
where
\begin{align*}
\epi (b) \coloneqq \{(t,x)\in [0,\infty]\times [-\infty ,\infty ] : x \geq b(t)\}.
\end{align*}
Since $b$ is lower semicontinuous $\epi (b)$ is closed in $[0,\infty]\times [-\infty ,\infty ]$. If $(X_t)_{t\geq 0}$ has right-continuous paths so has $(t,X_t)_{t\geq 0}$, and thus $\tau_b$ as a hitting time of a closed set is a stopping time, for instance see \cite{bass2018}.
\end{remark}

\begin{lemma}\label{cadlagandboundaryfunction}
Assume that $(X_t)_{t\geq 0}$ has right-continuous paths and $b$ is a boundary function. Then it holds
\begin{align*}
X_{\tau_b} \geq b(\tau_b)
\end{align*}
almost surely.
\end{lemma}
\begin{proof}
From the definition of the first-passage time we can find a (possibly random) sequence $s_n \searrow \tau_b$, such that $X_{s_n} \geq b(s_n)$ for all $n\in\N$. By the right-continuity it follows that
\begin{align*}
X_{\tau_b} = \lim_{n\to\infty} X_{s_n} \geq \liminf_{n\to\infty}b(s_n) \geq b(\tau_b),
\end{align*}
where the last inequality follows from the lower semicontinuity.
\end{proof}

For a boundary function $b$ define
\begin{align*}
\overline{\tau}_b \coloneqq \inf\{t\geq 0 : X_t \geq b(t)\}.
\end{align*}

\begin{lemma}\label{Anulova lemma 1}
Assume that $(X_t)_{t\geq 0}$ has right-continuous paths and is quasi-left-continuous. Furthermore, let $b_n \Hconv b$. Then on $\{\liminf_{n\to\infty} \tau_{b_n}>0\} \cup \{\tau_b = \overline{\tau}_b\}$ we have
\begin{align*}
\tau_b \leq \liminf_{n\to\infty} \tau_{b_n}
\end{align*}
almost surely.
\end{lemma}

\begin{proof}
We assume that $T\coloneqq \liminf_{n\to\infty} \tau_{b_n}<\infty$. Set $T_m \coloneqq \inf_{n \geq m}\tau_{b_n}$ and fix $m\in \N$. There is a sequence $(n_k)_{k\geq 1} \subseteq \{n\in \N: n\geq m\}$ (possibly random) such that $\lim_{k\to\infty}\tau_{b_{n_k}} = T_m$ and $\tau_{b_{n_{k}}} \geq \tau_{b_{n_{k+1}}}$ for all $k\in\N$. If $(n_k)_{k\in\N}$ is bounded, we can assume without loss of generality that $(n_k)_{k\geq 1}$ is a constant sequence and set $n^{(m)} \coloneqq n_1$. Then, from Lemma~\ref{cadlagandboundaryfunction}, we have $X_{T_m} \geq b_{n^{(m)}}(T_m)$ almost surely. On the other hand, if $n^{(m)} \coloneqq \lim_{k\to\infty} n_k =\infty$, due to the $\Gamma$-convergence and the right-continuity of the paths we have
\begin{align*}
b(T_m) \leq \liminf_{k\to\infty} b_{n_k}(\tau_{b_{n_k}}) \leq \liminf_{k\to\infty} X_{\tau_{b_{n_k}}}  = X_{T_m}
\end{align*}
almost surely. Note that the (possibly random) sequence of boundary functions
\begin{align*}
\tilde{b}_m \coloneqq \begin{cases} b_{n^{(m)}} &: n^{(m)} < \infty,\\ b &: n^{(m)} = \infty, \end{cases}
\end{align*}
$\Gamma$-converges to $b$ as $m\to\infty$ since $n^{(m)} \geq m$. Now, since $T = \lim_{m\to\infty} T_m$, we have, by the $\Gamma$-convergence and the quasi-left-continuity, that
\begin{align*}
b(T) \leq \liminf_{m\to\infty}  \tilde{b}_m (T_m) \leq \liminf_{m\to\infty} X_{T_m} = X_T
\end{align*}
almost surely. If $T>0$, it follows directly that $\tau_b \leq T = \liminf_{n\to\infty} \tau_{b_n}$. Generally, it follows $\overline{\tau}_b \leq T$, which concludes the proof.
\end{proof}

\begin{lemma}\label{liminf lemma}
Let $b_n \Hconv b$. Let $t\in (0,\infty)$ and assume that
\begin{align*}
\liminf_{s\searrow t}  b(s) = b(t).
\end{align*}
Then there exists a sequence $t_n \to t$ with $t_n > t$ such that
\begin{align*}
b_n(t_n) \to b(t).
\end{align*}
\end{lemma}
\begin{proof}
Since $\liminf_{s\searrow t} b(s) = b(t)$ there is a sequence $r_m \to t$ with $r_m > t$ such that $b(r_m) \to b(t)$ as $m\to\infty$. Since $b_n \Hconv b$ for every $m\in \N$ there is a sequence $r_n^m \to r_m$ such that $b_n (r_n^m) \to b(r_m)$. Without loss of generality we can assume that $r_n^m > t$. We now define two sequences $(m_k)_{k\in\N}$ and $(n_k)_{k\in\N}$ by a recursive scheme. For $k\in \N$ assume that $m_1, \ldots , m_{k-1}$ and $n_1, \ldots , n_{k-1}$ are already defined. Then let $m_k > m_{k-1}$ be large enough such that
\begin{align*}
\max (r_m -t , \abs{b(r_m) -b(t)}) \leq \frac{1}{k}\qquad \forall m \geq m_k.
\end{align*}
Further, let $n_k > n_{k-1}$ be large enough such that
\begin{align*}
\max( \abs{r_n^{m_k}-r_{m_k}} ,\abs{b_n (r_n^{m_k}) - b(r_{m_k})}) \leq \frac{1}{k}\qquad \forall n \geq n_k.
\end{align*}
Now define for $n\in\N$ the sequence
\begin{align*}
t_n \coloneqq \sum_{k=1}^\infty r_n^{m_k} \IF{n_k, \ldots , n_{k+1}-1}(n).
\end{align*}
Let $\varepsilon >0$. Choose $k\in\N$ such that $\frac{2}{k} < \varepsilon $. Let $n\in\N$. Then, if $n \in \{n_{\tilde{k}}, \ldots , n_{\tilde{k}+1}-1\}$ for $\tilde{k}\geq k$, we have
\begin{align*}
\abs{t_n-t} \leq \abs{r_n^{m_{\tilde k}} - r_{m_{\tilde k}}} +\abs{r_{m_{\tilde k}} - t}\leq \frac{1}{\tilde{k}} + \frac{1}{\tilde{k}} \leq \frac{2}{k} < \varepsilon.
\end{align*}
and
\begin{align*}
\abs{b_n(t_n) -b(t)} \leq \abs{b_n(r_n^{m_{\tilde{k}}}) -b(r_{m_{\tilde{k}}})} + \abs{b(r_{m_{\tilde{k}}})-b(t)} \leq \frac{1}{\tilde{k}} + \frac{1}{\tilde{k}} \leq \frac{2}{k} < \varepsilon.
\end{align*}
This shows that eventually that $t_n \to t$ with $t_n > t$ and
\begin{align*}
b_n(t_n) \to b(t)
\end{align*}
as $n\to\infty$.
\end{proof}

For a boundary function $b$ and $\varepsilon >0$ we interpret $b+\varepsilon$ as the boundary function given by $(b+\varepsilon) (t) = b(t) +\varepsilon$.

\begin{lemma}\label{Anulova lemma 2}
Assume that $(X_t)_{t\geq 0}$ has right-continuous paths and is quasi-left-continuous. Furthermore, let $b_n \Hconv b$ and $\varepsilon>0$. Then
\begin{align*}
\limsup_{n\to\infty} \tau_{b_n} \leq  \tau_{b+\varepsilon}
\end{align*}
almost surely.
\end{lemma}
\begin{proof}
According to Lemma~\ref{discontinuity lemma} the set
\begin{align*}
S_b \coloneqq \left\{ t\in [0,\infty) : \liminf_{s\searrow t} b(s) > b(t) \right\}
\end{align*}
is countable. By setting $b(0) \coloneqq \liminf_{s\searrow 0} b(s)$ (which does not affect $\tau_{b+\varepsilon}$) we can assume that $0\notin S_b$. Since $b$ is lower semicontinuous this means that for every $t\in [0,\infty) \setminus S_b$ we have $\liminf_{s\searrow t} b(s) = b(t)$. Since $(X_t)_{t\geq 0}$ is quasi-left-continuous and has right-continuous paths we have that
\begin{align}\label{eq:Anulova lemma 2}
X_{t-} = X_t \qquad \forall t\in S_b
\end{align}
almost surely. Assume without loss of generality that $\tau_{b+\varepsilon}<\infty$. By Lemma~\ref{cadlagandboundaryfunction} we have $X_{\tau_{b+\varepsilon}} \geq b(\tau_{b+\varepsilon}) +\varepsilon$. Due to the $\Gamma$-convergence we can choose a converging sequence $t_n \to \tau_{b+\varepsilon}$ (possibly random) such that $b_n (t_n) \to b(\tau_{b+\varepsilon})$ as $n\to\infty$. We distinguish two cases.

If $\tau_{b+\varepsilon} \in S_b$, by \eqref{eq:Anulova lemma 2}, we can assume that $X_{\tau_{b+\varepsilon}-} = X_{\tau_{b+\varepsilon}}$. We have therefore
\begin{align*}
\lim_{n\to\infty} X_{t_n} =X_{\tau_{b+\varepsilon}}.
\end{align*}

If $\tau_{b+\varepsilon} \in [0,\infty) \setminus S_b$ due to Lemma~\ref{liminf lemma}, we can assume that $t_n > \tau_{b+\varepsilon}$ for all $n\in\N$. Thus, since $(X_t)_{t\geq 0}$ has right-continuous paths we have
\begin{align*}
\lim_{n\to\infty} X_{t_n} = X_{\tau_{b+\varepsilon}}.
\end{align*}
Let $N\in\N$ (possibly random) be large enough such that for every $n\geq N$ we have
\begin{align*}
b_n (t_n ) \leq b(\tau_{b+\varepsilon}) + \frac{\varepsilon}{2} \quad \text{and} \quad X_{t_n } \geq b(\tau_{b+\varepsilon}) + \frac{\varepsilon}{2}.
\end{align*}
Since $t_n > 0$ it follows that $\tau_{b_n} \leq t_n$. In particular, we have
\begin{align*}
\limsup_{n\to\infty}\tau_{b_n} \leq \limsup_{n\to\infty} t_n =  \tau_{b+\varepsilon}.
\end{align*}
This shows the desired statement.
\end{proof}

\begin{remark}\label{remark Anulova lemma 3}
Let $b$ be a boundary function. Since $\tau_b' \geq \tau_b$, if for $(X_t)_{t\geq 0 }$ it holds that $\tau_b \overset{\text{d}}{=} \tau_b'$, then we have in fact $\tau_b = \tau_b'$ almost surely.
\end{remark}

\begin{lemma}\label{Anulova lemma 3}
Let $b$ be a boundary function. Then
\begin{align*}
\lim_{\varepsilon \searrow 0} \tau_{b+\varepsilon} = \tau_b'
\end{align*}
almost surely.
\end{lemma}
\begin{proof}
Note that $\tau_{b+\varepsilon}$ is decreasing in $\varepsilon$ and bounded from below by $\tau_b'$. Thus the following limit exists and fulfills
\begin{align*}
\lim_{\varepsilon \searrow 0 }\tau_{b+\varepsilon} \geq \tau_b'
\end{align*}
almost surely. Thus it is left to show that $\lim_{\varepsilon \searrow 0 }\tau_{b+\varepsilon} \leq \tau_b'$ almost surely. Without loss of generality we assume that $ \tau_b' < \infty$. On this event we have that there exists a sequence $t_n \searrow \tau_b'$ (possibly random) with
\begin{align*}
X_{t_n} > b(t_n) \quad \forall n \in \N.
\end{align*}
If $\tau_b'=0$ this sequence fulfills $t_n > 0$ for every $n\in\N$. Let $\delta > 0$. Then there exists $n\in \N$ (possibly random) such that $t_n \leq \tau_b' + \delta$. Set
\begin{align*}
\varepsilon_n \coloneqq \frac{1}{2} \left( X_{t_n} - b(t_n) \right) >0
\end{align*}
for which holds
\begin{align*}
X_{t_n} > b(t_n) + \varepsilon_n.
\end{align*}
Therefore $\tau_{b+\varepsilon_n}  \leq t_n \leq \tau_b' +\delta$. Thus we have that
\begin{align*}
\lim_{\varepsilon \searrow 0 }\tau_{b+\varepsilon}\leq \tau_{b+\varepsilon_n} \leq \tau_b' +\delta.
\end{align*}
Letting $\delta \searrow 0$ yields that $\lim_{\varepsilon \searrow 0 }\tau_{b+\varepsilon} \leq \tau_b' $. This proves the desired statement.
\end{proof}

\begin{lemma}\label{Anulova lemma 4}
Assume that $(X_t)_{t\geq 0}$ has right-continuous paths and is quasi-left-continuous. Let $b$ be a boundary function and assume $\tau_{b} \overset{\text{d}}{=} \tau_{b}'$. Further, let $b_n \Hconv b$. Then on $\{\liminf_{n\to\infty} \tau_{b_n}>0\} \cup \{\tau_b = \overline{\tau_b}\}$ we have
\begin{align*}
\lim_{n\to\infty}\tau_{b_n} =\tau_b
\end{align*}
almost surely.
\end{lemma}
\begin{proof}
By combining Lemma~\ref{Anulova lemma 1}, Lemma~\ref{Anulova lemma 2}, Lemma~\ref{Anulova lemma 3} and Remark~\ref{remark Anulova lemma 3} we have
\begin{align*}
\tau_b \leq \liminf_{n\to\infty} \tau_{b_n} \leq \limsup_{n\to\infty} \tau_{b_n} \leq \lim_{\varepsilon \searrow 0} \tau_{b+\varepsilon} = \tau_b
\end{align*}
almost surely on $\{\liminf_{n\to\infty} \tau_{b_n}>0\} \cup \{\tau_b = \overline{\tau_b}\}$, which concludes the proof.
\end{proof}

\begin{remark}
If $b_n \Hconv b$, in general it is not true that for every $s>0$ also $b_n|_s \Hconv b|_s$. The following Lemma~\ref{s truncation lemma} gives a sufficient condition on $b$ such that the convergence is preserved. By Lemma~\ref{discontinuity lemma} we will see that we can find arbitrarily small $s>0$ fulfilling this sufficient condition.
\end{remark}

\begin{lemma}\label{s truncation lemma}
Let $b_n \Hconv b$ and $s>0$. If $\liminf_{t\searrow s} b(t) = b(s)$, then $b_n|_s \Hconv b|_s$.
\end{lemma}
\begin{proof}
If $t\in [0,s)$ for every sequence $t_n\to t$ we have
\begin{align*}
\lim_{n\to\infty}b_n\vert_s (t_n) = b\vert_s(t) = \infty.
\end{align*}
If $t\in [s,\infty ]$ and $t_n\to t$, then
\begin{align*}
\liminf_{n\to\infty} b_n\vert_s (t_n) \geq \liminf_{n\to\infty} b_n (t_n) \geq b(t) = b\vert_s (t). 
\end{align*}
Furthermore, due to the $\Gamma$-convergence for $t\in [s,\infty]$ there is a sequence $t_n \to t$ such that $b_n(t_n) \to b(t)$. If $t=s$, by Lemma~\ref{liminf lemma} and our assumption, we can assume that $t_n \geq t=s$ for every $n$. Therefore we have
\begin{align*}
\lim_{n\to\infty} b_n\vert_s (t_n) = \lim_{n\to\infty} b_n (t_n) = b(t) = b\vert_s (t). 
\end{align*}
Hence it holds $b_n\vert_s \Hconv b\vert_s$.
\end{proof}

\begin{lemma}\label{Anulova lemma 5}
Let $b$ be a boundary function. Then
\begin{align*}
\tau_b = \lim_{s\searrow 0} \tau_{b|_s}
\end{align*}
almost surely.
\end{lemma}
\begin{proof}
The random variable $\tau_{b|_s}$ is monotone decreasing in $s$ and bounded from below by $\tau_b$. Thus the limit exists and we have $\tau_b \leq \lim_{s\searrow 0} \tau_{b|_s}$. Without loss of generality assume $\tau_b < \infty$. For $m\in \N$ there exists a time $t> 0$ (possibly random) such that $t\in [\tau_b , \tau_b + \frac{1}{m} )$ and $X_t \geq b(t)$. Thus, if $r \in (0 , t)$, then
\begin{align*}
X_t \geq b(t) = b|_{r} (t),
\end{align*} 
which means that $ \tau_{b|_r}  \leq t < \tau_b + \frac{1}{m}$. Consequently,
\begin{align*}
\tau_b \leq \lim_{s\searrow 0} \tau_{b|_s} \leq \tau_{b|_r} \leq \tau_b + \frac{1}{m}.
\end{align*}
By $m\to\infty$ we obtain $ \lim_{s\searrow 0} \tau_{b|_s} =  \tau_b$ almost surely.
\end{proof}

\begin{remark}\label{Anulova lemma 5 variant}
It is analogous to show that $\tau_b' = \lim_{s\searrow 0} \tau_{b\vert_s}'$ almost surely.
\end{remark}

\begin{proof}[Proof of Proposition~\ref{Anulova convergence}]
For $s>0$ and $\varepsilon >0$ we have that
\begin{align*}
&\PR{\abs{\tau_{b_n} - \tau_b}>\varepsilon}\\
&\leq \PR{ \abs{\tau_{b_n|_s} - \tau_{b_n} }>\frac{\varepsilon}{3}} + \PR{ \abs{\tau_{b_n|_s} - \tau_{b|_s} }>\frac{\varepsilon}{3}}+ \PR{ \abs{\tau_{b|_s} - \tau_{b} }>\frac{\varepsilon}{3}}\\
&\leq \PR{\tau_{b_n} \leq s} + \PR{ \abs{\tau_{b_n|_s} - \tau_{b|_s} }>\frac{\varepsilon}{3}}+ \PR{ \abs{\tau_{b|_s} - \tau_{b} }>\frac{\varepsilon}{3}},
\end{align*}
where we have used that
\begin{align*}
\PR{ \abs{\tau_{b_n|_s} - \tau_{b_n} }>\frac{\varepsilon}{3}} \leq \PR{ \tau_{b_n|_s} \neq \tau_{b_n} }\leq  \PR{\tau_{b_n} \leq s}.
\end{align*}
Due to Lemma~\ref{discontinuity lemma} we can choose arbitrarily small $s>0$ such that $\liminf_{t\searrow s}b(t) = b(s)$. Lemma~\ref{s truncation lemma} shows that $b_n\vert_s \Hconv b\vert_s$. Moreover, we have almost surely
\begin{align*}
\liminf_{n\to\infty} \tau_{b_n|_s} \geq s >0
\end{align*}
and by assumption $\tau_{b\vert_s} \overset{\text{d}}{=} \tau_{b\vert_s}'$. Hence we can apply Lemma~\ref{Anulova lemma 4} and obtain
\begin{align*}
\lim_{n\to\infty} \tau_{b_n|_s} = \tau_{b|_s}
\end{align*}
almost surely. By Lemma~\ref{Anulova lemma 5} we have
\begin{align*}
\lim_{s\searrow 0}\tau_{b|_s}  = \tau_b.
\end{align*}
almost surely. Thus we have that
\begin{align*}
&\limsup_{n\to\infty}\PR{\abs{\tau_{b_n} - \tau_b}>\varepsilon}\\
&\leq \lim_{s\searrow 0} \left( \limsup_{n\to\infty} \left(\PR{\tau_{b_n} \leq s} + \PR{ \abs{\tau_{b_n|_s} - \tau_{b|_s} }>\frac{\varepsilon}{3}}\right) +  \PR{ \abs{\tau_{b|_s} - \tau_{b} }>\frac{\varepsilon}{3}}\right)
 \\
&= \lim_{s\searrow 0}\limsup_{n\to\infty} \PR{\tau_{b_n} \leq s}=  0,
\end{align*}
where the last equality comes from our assumption. This yields the statement.
\end{proof}

\begin{remark}\label{compactness}
Analogously as in Theorem~3.1 of \cite{kk2022} it can  be shown that the $\Gamma$-convergence coincides with the convergence of the epigraphs in the Hausdorff metric. From Proposition~2.1.3 of \cite{k2022} follows that every sequence of boundary functions has a convergent subsequence.
\end{remark}

\begin{proof}[Proof of Theorem~\ref{existence}]
For $n\in\N$ let $t_k^n \coloneqq k 2^{-n}$ with $k\in \N_0$. We will inductively define a boundary function $b_n$, which has only finite values at the discrete timepoints $t_k^n$. For $k\in \N$ let us assume $b_n(t_1^n), \ldots , b_n (t_{k-1}^n)$ are already defined. Since $ \PR{X_t \in \cdot\;}$ is diffuse, we can choose a value $b_n (t_k^n) \in [-\infty , \infty]$ such that
\begin{align*}
\PR{X_{t_k^n} < b_n (t_k^n) , \ldots , X_{t_1^n} < b_n (t_1^n)} = \PR{\xi > t_k^n} 
\end{align*}
with $b_n (t_k^n) = -\infty$ if $\PR{\xi > t_k^n} =0$. By setting $b_n (t) = \infty$ for all $t \notin \{k2^{-n} : k\in \N\}$, we obtain a lower semicontinuous function $b_n$. Note, that then by definition
\begin{align*}
\PR{\tau_{b_n} > t} = \PR{\xi > \lfloor t2^n \rfloor 2^{-n}}, \qquad \forall t\geq 0.
\end{align*}
This implies that
\begin{align*}
\tau_{b_n} \overset{\text{d}}{\rightarrow} \xi
\end{align*}
as $n\to\infty$. By the compactness of the set of boundary functions, see Remark~\ref{compactness}, there is a lower semicontinuous function $b$, and a subsequence $N\subset \N$ such that
\begin{align*}
b_n \Hconv b
\end{align*}
along $n\in N$. By assumption we have $\tau_{b\vert_s}\overset{\text{d}}{=} \tau_{b\vert_s}'$ for every $s>0$. Moreover, by assumption we have that $(X_t)_{t\geq 0}$ has right-continuous paths and is quasi-left-continuous. From Proposition~\ref{Anulova convergence} and Remark~\ref{Anulova convergence remark} we obtain that $\tau_{b_n} \overset{\Pr}{\to} \tau_b$ in probability along $n\in N$. This implies that $\tau_{b} \overset{\text{d}}{=} \xi$.
\end{proof}

\section{Uniqueness: Proof of Theorem~\ref{uniqueness}}\label{s:uniqueness}

Let us explain beforehand the role of the conditions in the proof of Theorem~\ref{uniqueness}. Conditions \ref{uniqueness: markov assumption} and \ref{uniqueness: order preservation assumption} allow to construct a boundary function which is a lower bound for any other solution. Conditions \ref{existence: diffusitivy assumption}, \ref{existence: continuity assumption}, \ref{existence: first-passage time assumption} will yield that this lower bound is a solution. Condition \ref{uniqueness: support assumption} will allow to infer that this lower bound is the unique solution.

Let $E\subseteq \R$ be an interval with $\overline E = [L,R]$ and $L,R \in [-\infty ,  \infty ]$. Assume that $ \PR{X_t \in \cdot\;}$ is diffuse for any $t>0$ and assume that $\PR{\tau_R < \infty} = 0$.

Let $(t_k^n)_{n\in\N, k\in \{0,1, \ldots , m_n\}} \subset [0,\infty )$ with $m_n \in \N\cup \{\infty\}$ be such that
\begin{align*}
0 = t_0^n < t_1^n < \ldots t_{m_n}^n.
\end{align*}
For $k\in \N$, such that $\PR{\xi > t_k^n} >0$ suppose $q_1^n , \ldots q_{k-1}^n$ are already defined. Since $\PR{X_{t_k^n} \in \cdot\;}$ is diffuse and $\PR{X_{t_k^n} < R}=1$ we can choose $q_k^n \in \overline{E}$ such that
\begin{align*}
\PR{X_{t_k^n} < q_k^n , X_{t_{k-1}^n} < q_{k-1}^n, \ldots , X_{t_1^n} < q_1^n} = \PR{\xi > t_k^n}
\end{align*}
and
\begin{align*}
\PR{X_{t_k^n} < q , X_{t_{k-1}^n} < q_{k-1}^n, \ldots , X_{t_1^n} < q_1^n} < \PR{\xi > t_k^n}
\end{align*}
for any $q < q_k^n$. Note that we have $q_k^n > \inf E = L$ since $\PR{\xi > t_k^n} >0$. For $k\in \N$ with $\PR{\xi > t_k^n}=0$ we set $q_k^n \coloneqq \inf E =L$. By setting
\begin{align}\label{lowerbarrierapprox}
b_n (t) \coloneqq \begin{cases}
q_k^n &: t = t_k^n,\\
\sup E = R &: t \notin \{t_k^n : k\in \N\},
\end{cases}
\end{align}
we obtain a boundary function $b_n$ with values in $\overline{E}$. Note that by the definition of $q_k^n$ we obtain that
\begin{align*}
b_n(t_k^n) = \sup \supp ( \PRb{X_{t_k^n} \in\cdot}{\tau_{b_n}>t_k^n} ).
\end{align*}

\begin{remark}
For Brownian motion on $\R$ this discretization appeared in \cite{Anu80} and \cite{kk2022}, and implicitely in \cite{Eks16}, \cite{demasi2019}. In \cite[Lemma~4.1]{kk2022} and \cite[Theorem~5]{demasi2019} it led to statements which are special cases of the following Lemma~\ref{lower barrier approximation lemma}.
\end{remark}

\begin{lemma}\label{lower barrier approximation lemma}
Let $E\subseteq \R$ be an interval with $\overline E = [L,R]$. Let $\mu$ be a probability measure on $E$. Assume \ref{uniqueness: markov assumption}, \ref{uniqueness: order preservation assumption} and \ref{existence: diffusitivy assumption} and $\PR{\tau_R < \infty} = 0$ with $\Pr \coloneqq \Pr_\mu$. Let $b$ be a boundary function with values in $\overline{E}$ such that $\tau_b  \overset{\text{d}}{=} \xi$. Then for fixed $n\in \N$ we have
\begin{align*}
\PRb{X_{t_k^n} \in\cdot}{\tau_{b_n}>t_k^n} \st \PRb{X_{t_k^n} \in\cdot}{\tau_{b}>t_k^n} \qquad \forall k\in\N : \PR{\xi > t_k^n} >0.
\end{align*}
In particular, it follows that $b_n (t_k^n) \leq b(t_k^n)$ for all $k\in \N$.
\end{lemma}

In order to prove Lemma~\ref{lower barrier approximation lemma} we need one more tool. For a probability measure $\mu$ and $\alpha \in (0,1]$ define for $A\subseteq \R$ measurable
\begin{align*}
T_\alpha (\mu)(A) \coloneqq \frac{\mu ( A \cap ( -\infty  , q_\alpha (\mu) ] )}{\mu ((-\infty , q_\alpha (\mu )])},
\end{align*}
where
\begin{align*}
q_\alpha (\mu) \coloneqq \inf\{c \in \R : \mu ( (-\infty , c]) \geq \alpha\}.
\end{align*}

The following statement is the one-sided version of Lemma~3.3 in \cite{kk2022}. In the presented generality we will use it in Section~\ref{section:comparison principle}.

\begin{lemma}\label{truncation preserves order}
Let $\mu$, $\nu$ be probability measures and such that $\mu$ is diffuse. Then for $\alpha_1 , \alpha_2 \in (0,1]$ with $\alpha_1 \leq \alpha_2$ we have that $\mu \st \nu$ implies $T_{\alpha_1} (\mu) \st T_{\alpha_2} (\nu)$.
\end{lemma}
\begin{proof}
Since $\mu$ is diffuse we have that $\mu ((-\infty, q_{\alpha_1} (\mu)])=\alpha_1$. Assume that $\mu \st \nu$. Then by the definitions we have $q_{\alpha_1} (\mu) \leq q_{\alpha_2} (\nu)$ and it suffices to consider the case $c\leq q_{\alpha_1} (\mu)$. Since $\nu ((-\infty, q_{\alpha_2} (\nu)]) \geq \alpha_2$, we have
\begin{align*}
T_{\alpha_1} (\mu ) ((-\infty , c]) &= \frac{\mu ((-\infty , c])}{\mu ((-\infty,q_{\alpha_1} (\mu)])} = \frac{\mu ((-\infty , c])}{\alpha_1} \geq \frac{\mu ((-\infty , c])}{\alpha_2} \\
&\geq \frac{\mu ((-\infty , c])}{\nu ((-\infty, q_{\alpha_2} (\nu)])} \geq \frac{\nu ((-\infty , c])}{\nu ((-\infty, q_{\alpha_2} (\nu)])} = T_{\alpha_2} (\nu) ((-\infty, c]). 
\end{align*}
This shows $T_\alpha (\mu) \st T_\alpha (\nu)$.
\end{proof}

For a probability measure $\mu$ we introduce the mapping $P_t$ by
\begin{equation}\label{transition kernel}
P_t (\mu) \coloneqq \PRx{\mu}{X_t \in \cdot \;}.
\end{equation}

\begin{proof}[Proof of Lemma~\ref{lower barrier approximation lemma}]
Essentially, we follow the lines of the proof of Lemma~4.1 in \cite{kk2022}, which was conducted in the case of reflected Brownian motion. Fix $n\in\N$. We abbreviate
\begin{align*}
\mu_k^{n} \coloneqq \PRb{X_{t_k^n} \in\cdot}{\tau_{b_n}>t_k^n}, \quad \mu_k \coloneqq \PRb{X_{t_k^n} \in\cdot}{\tau_{b}>t_k^n}.
\end{align*}
Note that from the Markov property it follows that
\begin{align}\label{lower barrier approximation lemma proof 1}
\mu_k^{n} = T_{\alpha_k^n} \circ P_{t_k^n -t_{k-1}^n} \circ \ldots T_{\alpha_1^n} \circ P_{t_1^n} (\mu),
\end{align}
where
\begin{align*}
\alpha_k^n = \frac{\PR{\xi > t_k^n}}{\PR{\xi > t_{k-1}^n}}.
\end{align*}
We will prove the statement by induction over $k\in \N_0$ with $\PR{\xi > t_k^n} >0$, by comparing the mappings
\begin{align*}
H_k^n (\nu ) \coloneqq T_{\alpha_k^n} \circ P_{t_k^n -t_{k-1}^n} (\nu),
\end{align*}
where $\nu$ is a probability measure on $E$, and
\begin{align*}
H_k (\nu ) \coloneqq \PRbx{\nu}{X_{t_k^n -t_{k-1}^n} \in \cdot}{\tau_{b^{t_{k-1}^n}} > t_k^n -t_{k-1}^n},
\end{align*}
where we used the notation $b^s (t) \coloneqq b (t+s)$ for $s>0$. It follows by \eqref{lower barrier approximation lemma proof 1} and the Markov property that
\begin{align*}
H_k^n (\mu_{k-1}^n) = \mu_k^n \qquad \text{and} \qquad H_k (\mu_{k-1}) = \mu_k.
\end{align*}
We now claim that we have
\begin{align}\label{lower barrier approximation lemma proof 2}
H_k^n (\mu_{k-1}) \st H_k (\mu_{k-1}).
\end{align}
Using the Markov property we obtain
\begin{align*}
P_{t_k^n -t_{k-1}^n} (\mu_{k-1}) = \PRbx{\mu}{X_{t_k^n}\in\cdot \;}{\tau_b > t_{k-1}^n }.
\end{align*}
This shows that $P_{t_k^n -t_{k-1}^n} (\mu_{k-1})$ is diffuse and we have
\begin{align*}
P_{t_k^n -t_{k-1}^n} (\mu_{k-1})((-\infty, q_{\alpha_k^n} (P_{t_k^n -t_{k-1}^n} (\mu_{k-1}))]) = \alpha_k^n
\end{align*}
and, by the Markov property and the fact that $\tau_b \overset{\text d} = \xi$, we have
\begin{align*}
\PRx{\mu_{k-1}}{ \tau_{b^{t_{k-1}^n}} > t_k^n -t_{k-1}^n} = \alpha_k^n.
\end{align*}
Therefore, if $c\leq q_{\alpha_k^n} (P_{t_k^n -t_{k-1}^n} (\mu_{k-1}))$, we have
\begin{align*}
&H_k^n (\mu_{k-1}) ((-\infty ,c ]) = \frac{P_{t_k^n -t_{k-1}^n} (\mu_{k-1}) ((-\infty , c])}{\alpha_k^n}\\
& = \frac{\PRx{\mu_{k-1}}{X_{t_k^n -t_{k-1}^n}  \leq c}}{\alpha_k^n} \geq \frac{\PRx{\mu_{k-1}}{X_{t_k^n -t_{k-1}^n}  \leq c, \tau_{b^{t_{k-1}^n}} > t_k^n -t_{k-1}^n}}{\alpha_k^n}\\
&= \frac{\PRx{\mu_{k-1}}{X_{t_k^n -t_{k-1}^n}  \leq c, \tau_{b^{t_{k-1}^n}} > t_k^n -t_{k-1}^n}}{\PRx{\mu_{k-1}}{ \tau_{b^{t_{k-1}^n}} > t_k^n -t_{k-1}^n}}\\
&= \PRbx{\mu_{k-1}}{X_{t_k^n -t_{k-1}^n}  \leq c}{ \tau_{b^{t_{k-1}^n}} > t_k^n -t_{k-1}^n} = H_k (\mu_{k-1}) ((-\infty , c]).
\end{align*}
This shows the claim. Now let us assume that $\mu_{k-1}^n \st \mu_{k-1}$. Then Lemma~\ref{truncation preserves order} the fact that $P_t$ preserves the usual stochastic order and the claim from \eqref{lower barrier approximation lemma proof 2} yield
\begin{align*}
\mu_k^n = H_k^n (\mu_{k-1}^n) \st  H_k^n (\mu_{k-1}) \st H_k (\mu_{k-1}) = \mu_k.
\end{align*}
Since $\mu_0^n = \mu_0$ the desired ordering follows by induction. From the ordering it follows that $\mu_k^n ( (-\infty , c]) \geq  \mu_k ((-\infty , c])$ for all $c\in \R$, hence
\begin{align*}
b_n (t_k^n) = \sup \supp (\mu_k^n) \leq \sup \supp (\mu_k) \leq b(t_k^n)
\end{align*}
for $k\in\N$ with $\PR{\xi > t_k^n}>0$. Since $b_n (t_k^n) = \inf E$ for $k\in \N$ with $\PR{\xi > t_k^n}=0$, the proof is finished.
\end{proof}

Let $b$ be a boundary function with values in $\overline{E}$. We now introduce a discretization technique for $b$, which was already used in \cite{Eks16} and \cite{chen2011} for the case of Brownian motion. We use an adapted version. Let $D(b)$ be an arbitrary countable set and $D_n(b) \subset D_{n+1} (b)$ finite, such that $\bigcup_{n\in \N} D_n (b) = D(b)$. For $n\in \N$ define for every $k\in \N$
\begin{align*}
\tilde t_k^n \coloneqq \inf \left\{ t \in [k2^{-n} , (k+1)2^{-n} ] : b(t) = \inf_{s\in [k2^{-n} , (k+1)2^{-n} ]} b(s) \right\}.
\end{align*}
Set $A^1_n (b) \coloneqq \{\tilde t_k^n : k\in \{1, 2 , \ldots , n2^n\}\}\}$. Furthermore, let $(s_n)_{n\in\N}$ be an enumeration of $\{s\in [0,\infty ) : \PR{\tau_b = s} > 0\}$. Set $A_n^2 (b) \coloneqq \{s_1 , \ldots , s_n\}$. Finally, set
\begin{align*}
A_n (b) \coloneqq A_n^1 (b) \cup A_n^2 (b) \cup D_n (b).
\end{align*}
By choosing $D_n (b) = \emptyset$ we end up with the construction used in \cite{Eks16}. Note that by $A_n^1 (b) \subset A_{n+1}^1 (b)$ we have
\begin{align*}
A_n (b) \subset A_{n+1} (b).
\end{align*}
For $n\in\N$ let us define the boundary function
\begin{align}\label{eksapprox}
\hat b_n (t) \coloneqq \begin{cases} b(t) &: t\in A_n (b),\\ R &: t\notin A_n (b). \end{cases}
\end{align}

\begin{lemma}\label{eksapprox gammaconverges}
For a boundary function $b$ it holds $\hat b_n \Hconv b$ as $n\to\infty$.
\end{lemma}
\begin{proof}
Let $t\in [0,\infty]$. Assume $t_n \to t$. Then
\begin{align*}
\hat b_n(t_n) \geq b(t_n ).
\end{align*}
Thus,
\begin{align*}
\liminf_{n\to\infty } \hat b_n (t_n) \geq \liminf_{n\to\infty} b (t_n) \geq b(t).
\end{align*}
For the second part of the $\Gamma$-convergence we distinguish two cases. Let us first assume that $t\in \bigcup_{n\in\N} A_n (b)$. Then for $N$ large enough we have $t\in A_n(b)$ for all $n\geq N$. Hence,
\begin{align*}
\lim_{n\to\infty} \hat b_n (t) = \lim_{n\to\infty} b(t) = b(t).
\end{align*}
Assume that $t\notin A_n (b)$ for all $n\in\N$. Let $n\in \N$ be large enough and $k_n (t) \in \N$ such that
\begin{align*}
k_n (t) 2^{-n} \leq t \leq (k_n(t) +1) 2^{-n}.
\end{align*}
Then we have
\begin{align*}
b(t) \geq \inf_{s \in [k_n(t) 2^{-n}, (k_n(t)+1)2^{-n}]} b(s) = b( \tilde t_{k_n(t)}^n).
\end{align*}
Now define $t_n \coloneqq \tilde t_{k_n (t)}^n$. We have $t_n \in A_n (b)$. It follows $t_n \to t$ and
\begin{align*}
b(t) &\geq \limsup_{n\to\infty} b (t_n) = \limsup_{n\to\infty} \hat b_n (t_n)\\
&\geq \liminf_{n\to\infty} \hat b_n (t_n) = \liminf_{n\to\infty} b(t_n) \geq b(t).
\end{align*}
This means that $\hat b_n (t_n) \to b(t)$. Altogether we obtain $\hat b_n \Hconv b$.
\end{proof}

\begin{lemma}\label{uniqueness for ordered}
Let $b_1$ and $b_2$ be boundary functions with values in $\overline{E}$ such that $\tau_{b_1}  \overset{\text{d}}{=} \tau_{b_2}  \overset{\text{d}}{\eqqcolon}\xi$ and $b_1 \leq b_2$. Let $t\in (0,t^\xi)$ and assume
\begin{align*}
b_2(t) = \sup \supp (\PRb{X_t \in \cdot}{\tau_{b_2} > t}).
\end{align*}
Then $b_1 (t) = b_2 (t)$.
\end{lemma}
\begin{proof}
Assume $b_1(t) < b_2(t)$. Then by the assumption for the support we would have
\begin{align*}
\PR{\xi > t} &= \PR{\tau_{b_2}> t} =  \PR{\tau_{b_2}> t , X_t < b_2 (t)}\\
&> \PR{\tau_{b_2}> t , X_t < b_1 (t)} \geq \PR{\tau_{b_1}> t, X_t < b_1(t)}\\
&= \PR{\tau_{b_1} > t} = \PR{\xi > t}.
\end{align*}
This contradiction shows $b_1(t) = b_2(t)$.
\end{proof}

\begin{lemma}\label{gammaconvergence preserves order}
Let $b^1$ and $b^2$ be boundary functions and $b_n^1 \Hconv b^1$ and $b_n^2 \Hconv b^2$. Assume that $b_n^1 \leq b_n^2$. Then $b^1 \leq b^2$.
\end{lemma}
\begin{proof}
Let $t\in [0,\infty]$. Let $t_n \to t$ such that $b_n^2(t_n) \to b^2(t)$. Then
\begin{align*}
b^2(t) = \lim_{n\to\infty} b^2 (t_n ) \geq \liminf_{n\to\infty} b^1 (t_n) \geq b^1 (t).
\end{align*}
This finishes the proof.
\end{proof}

\begin{proof}[Proof of Theorem~\ref{uniqueness}]
Let $b, \beta$ be boundary functions with values in $\overline{E}$ such that $\tau_b  \overset{\text{d}}{=} \tau_\beta \overset{\text{d}}{=} \xi$. Recall the construction from \eqref{eksapprox} for $b$ and $\beta$, respectively. In particular, recall $A^1_n (b), A^2_n (b)$ and $A^1_n (\beta) , A^2_n (\beta)$. We can choose
\begin{align*}
D_n (b) \coloneqq A_n^1 (\beta ) \cup A_n^2 (\beta ), \qquad D_n (\beta ) \coloneqq A^1_n (b) \cup A^2_n (b).
\end{align*}
This means that in the construction from \eqref{eksapprox} we have
\begin{align*}
A_n \coloneqq A_n (b) = A_n (\beta).
\end{align*}
Now since $\PR{X_t \in \cdot\;}$ is diffuse we can construct the boundary function $b_n$ from the construction \eqref{lowerbarrierapprox}, where we choose
\begin{align*}
\{t_0^n , t_1^n , \ldots , t^n_{m_n}\} = \{0\} \cup A_n.
\end{align*}
Due to the Markov property and the stochastic order preservation we can apply Lemma~\ref{lower barrier approximation lemma} for the solutions $b$ and $\beta$ separately but with the same set $A_n$ of discrete timesteps. With recalling \eqref{eksapprox} this leads to
\begin{align*}
b_n (t) \leq b(t) = \hat{b}_n (t) \quad \text{and} \quad b_n (t) \leq b(t) = \hat{\beta}_n (t) \qquad\forall t \in A_n.
\end{align*}
This means that $b_n \leq \hat{b}_n$ and $b_n \leq \hat{\beta}$ altogether. Now note that, by Remark~\ref{compactness}, there is a boundary function $b^+$ and a subsequence $N\subset \N$ such that
\begin{align*}
\lim_{n\in N} b_n \Hconv b^+.
\end{align*}
By Lemma~\ref{eksapprox gammaconverges} we have that
\begin{align*}
\hat{b}_n  \Hconv b \quad\text{and}\quad \hat{\beta}_n \Hconv \beta
\end{align*}
Thus we have by Lemma~\ref{gammaconvergence preserves order} that
\begin{align*}
b^+  \leq b \quad\text{and}\quad  b^+ \leq \beta.
\end{align*}
By the definition of $b_n$ from \eqref{lowerbarrierapprox} under $\Pr$ we have on the one hand that
\begin{align*}
\tau_{b_n} \overset{\text{d}}{\rightarrow} \xi
\end{align*}
as $n\to\infty$. By assumption we have $\tau_{b^+\vert s}= \tau_{b^+\vert s}'$ $\Pr$-a.s. for $s>0$ and $(X_t)_{t\geq 0}$ has $\Pr$-a.s. right-continuous paths and is quasi-left-continuous. By Proposition~\ref{Anulova convergence} and Remark~\ref{Anulova convergence remark} we obtain
\begin{align*}
\tau_{b_n} \overset{\Pr}\to \tau_{b^+}
\end{align*}
in probability. This means that $\tau_{b^+} \overset{\text{d}}{=} \xi$ under $\Pr$. Since $b_n$ has values in $\overline{E}$, by the definition of the $\Gamma$-convergence, it follows that $b^+$ is a boundary function with values in $\overline{E}$. By assumption we have
\begin{align*}
b(t) = \sup \supp \left(\PRb{X_t \in \cdot}{\tau_{b} > t}\right), \beta (t) = \sup \supp \left(\PRb{X_t \in \cdot}{\tau_{\beta} > t}\right)
\end{align*}
for $t \in I^\xi$. Hence Lemma~\ref{uniqueness for ordered} yields that
\begin{align*}
b(t) = b^+(t) = \beta (t)
\end{align*}
for every $t\in I^\xi$.
\end{proof}

\section{Comparison principle: Proof of Theorem~\ref{comparison principle}}\label{section:comparison principle}

The following proof essentially follows the lines of the proof of Theorem~2.2 of \cite{kk2022}. Due to its brevity we include it for completeness.

\begin{proof}[Proof of Theorem~\ref{comparison principle}]
Let $t_k^n \coloneqq k2^{-n}$ with $k\in\N_0$. For the measure $\Pr_{\mu_i}$ and the random variable $\xi_i$ let $b^i_n$ be the sequence of boundary function constructed in \eqref{lowerbarrierapprox}. For $k\in\N$ with $\PR{\xi_i > t_k^n}>0$ let
\begin{align*}
\alpha_k^{i,n} \coloneqq \frac{\PR{\xi_i > t_k^n}}{\PR{\xi_i > t_{k-1}^n}}.
\end{align*}
Since $\xi_1 \hr \xi_2$ we have that $\alpha_k^{1,n} \leq \alpha_k^{2,n}$. Recall the mapping $P_t (\mu ) = \PRx{\mu}{X_t \in \cdot \;}$ from \eqref{transition kernel}. Since $P_t$ preserves the order $\st$ and $\PRx{\mu_i}{X_t \in \cdot\;}$ are diffuse measures, we obtain by Lemma~\ref{truncation preserves order} that
\begin{align*}
&\PRbx{\mu_1}{X_{t_k^n} \in\cdot}{\tau_{b_n^1}>t_k^n} = T_{\alpha_k^{1,n}} \circ P_{t_k^n -t_{k-1}^n} \circ \ldots T_{\alpha_1^{1,n}} \circ P_{t_1^n} (\mu_1)\\
&\st T_{\alpha_k^{2,n}} \circ P_{t_k^n -t_{k-1}^n} \circ \ldots T_{\alpha_1^{2,n}} \circ P_{t_1^n} (\mu_2) = \PRbx{\mu_2}{X_{t_k^n} \in\cdot}{\tau_{b_n^2}>t_k^n}.
\end{align*}
This implies
\begin{align*}
b_n^1 (t_k^n) &= \sup \supp ( \PRbx{\mu_1}{X_{t_k^n} \in\cdot}{\tau_{b_n^1}>t_k^n} ) \\
&\leq \sup \supp ( \PRbx{\mu_2}{X_{t_k^n} \in\cdot}{\tau_{b_n^2}>t_k^n} ) = b_n^2 (t_k^n)
\end{align*}
for $k\in\N$ with $\PR{\xi_1 > t_k^n}>0$. Since $t^{\xi_1} \leq t^{\xi_2}$ this means that $b_n^1 \leq b_n^2$. Now let $b_1,b_2$ be accumulation points of the sequences $(b_n^1)_{n\in\N}$ and $(b_n^2)_{n\in\N}$ such that $N\subseteq \N$ is a subsequence with $b_n^i \Hconv b_i$ along $N$. Lemma~\ref{gammaconvergence preserves order} implies that
\begin{align*}
b_1 \leq b_2.
\end{align*}
As in the proof of Theorem~\ref{uniqueness} we have that Proposition~\ref{Anulova convergence} implies that $\tau_{b_i} \overset{\text{d}}{=} \xi_i$ under $\Pr_{\mu_i}$.
\end{proof}

\section{Conditions for Lévy processes: Proof of Theorem~\ref{Levy processes result}}\label{conditions for Levy}

In this section let $(X_t)_{t\geq 0}$ be a Lévy process on $\R$, where we allow $\PR{X_0 \in \cdot\; }$ to be an arbitrary probability measure on $\R$. We will show that under suitable conditions $(X_t)_{t\geq 0}$ fulfills the conditions of Theorem~\ref{uniqueness}, under which we established existence and uniqueness for the inverse first-passage time problem. This leads to the proof of Theorem~\ref{Levy processes result}, which is to be found at the end of the section. At first, we will collect the essential steps in preliminary statements. We begin with the fact that \ref{existence: diffusitivy assumption} already implies \ref{existence: first-passage time assumption} for Lévy processes.

\begin{proposition}\label{eks first passage time variant}
Let $(X_t)_{t\geq 0}$ be a Lévy process such that $\PR{X_1\in \cdot\;}$ is diffuse. Let $b$ be a boundary function. Then
\begin{align*}
\tau_b = \tau_b'
\end{align*}
almost surely.
\end{proposition}
The key idea for the proof of the statement is taken from Lemma~6.2 of \cite{Eks16}, where the statement was proved for Brownian motion in an a very similiar manner. For diffusions on $\R$ a corresponding statement was shown in \cite{chen2011}.
\begin{proof}
As first step we will assume that $b =b\vert_s$ for some $s>0$. Since $X_s$ is independent from the future increments and its law is diffuse, we have for $t\geq s$ that the law of
\begin{align*}
Z_t \coloneqq \sup_{r\in [s,t]} (X_r -b(r)) = X_s +\sup_{r\in [s,t]} (X_r - X_s - b(r))
\end{align*}
is diffuse. In particular it holds $\PR{Z_t=0}=0$. Recall that for a process with right-continuous paths we have $X_{\tau_b} \geq b(\tau_b)$ almost surely. Moreover, since $\inf_{t\in [0,s)} b(t) = \infty$ it holds that $\tau_b \geq s$. For $t\geq s$ we have
\begin{align*}
\{\tau_b \leq t\} \subseteq \{Z_t \geq 0\}, \quad \{Z_t > 0\} \subseteq \{\tau_b' \leq t\}.
\end{align*}
Consequently, we have for $t\geq s$ that
\begin{align*}
\PR{\tau_b \leq t } &= \PR{\tau_b \leq t , Z_t \geq 0} = \PR{\tau_b \leq t, Z_t > 0}\\
&= \PR{\tau_b \leq t, Z_t > 0, \tau_b' \leq t } =\PR{\tau_b \leq t,Z_t \geq 0, \tau_b' \leq t }\\
&= \PR{\tau_b \leq t, \tau_b' \leq t } = \PR{\tau_b'\leq t}.
\end{align*}
This shows $\tau_b \overset{\text{d}}= \tau_b'$.

For an arbitrary boundary function $b$ and $s>0$, since $\tau_{b\vert_s} \overset{\text{d}}= \tau_{b\vert_s}'$ and $\tau_{b\vert_s} \leq \tau_{b\vert_s}'$, it follows that $\tau_{b\vert_s} = \tau_{b\vert_s}'$ almost surely. Lemma~\ref{Anulova lemma 5} and Remark~\ref{Anulova lemma 5 variant} yield
\begin{align*}
\tau_b = \lim_{s\searrow 0}\tau_{b\vert_s} = \lim_{s\searrow 0} \tau_{b\vert_s}' =\tau_b'
\end{align*}
almost surely. This finishes the proof.
\end{proof}

For a measure $\mu$ on $\R$ we denote with $\suppex(\mu)$ the closure of $\supp (\mu)$ in $[-\infty , \infty]$. Recall the definition of a characteristic triple of a Lévy process in \eqref{characteristic triple}.

\begin{proposition}\label{levy process supportbedingung}
Let $(X_t)_{t\geq 0}$ be a Lévy process with characteristic triple $(a,\sigma^2 , \Pi )$ and $X_0 = 0$.
\begin{enumerate}
\item[(I)] If $0\in \supp (\Pi)$ and $(X_t)_{t\geq 0}$ is a subordinator without drift, then for every boundary function $b: [0,\infty ] \to [0,\infty ]$ and $t>0$ with $\PR{\tau_b > t} >0$ we have
\begin{align*}
\suppex (\PRb{X_t \in \cdot \;}{\tau_b > t}) = [0,b(t)].
\end{align*}
\item[(II)] If $0\in \supp (\Pi)$ and $(-X_t)_{t\geq 0}$ is a subordinator without drift, then for every boundary function $b: [0,\infty ] \to [-\infty ,0 ]$ with $\PR{\tau_b >0}>0$ and $t>0$ with $b(t) \leq b(u) $ for all $u\in [0,t]$ we have
\begin{align*}
\supp (\PRb{X_t \in \cdot \;}{\tau_b > t}) = (-\infty,b(t)].
\end{align*}
\item[(III)] If one of the following holds,
\begin{itemize}
\item[(i)] $(X_t)_{t\geq 0}$ has unbounded variation, i.e. $\sigma \neq 0$ or $\int_\R (1\wedge \abs{x}) \Pi (\intdiff x) = \infty$,
\item[(ii)] $0\in \supp (\Pi)$ and $\Pi ((-\infty , 0))>0$ and $\Pi ( (0 , \infty ))>0$,
\end{itemize}
then for every boundary function $b: [0,\infty ] \to [-\infty ,\infty ]$ and $t>0$ with $\PR{\tau_b > t} >0$ we have
\begin{align*}
\suppex (\PRb{X_t \in \cdot\;}{\tau_b>t} ) = (-\infty,b(t)].
\end{align*}
\end{enumerate}
\end{proposition}

The idea of the proof for Proposition~\ref{levy process supportbedingung} is to use the Lévy-It\^o decomposition and extract suitable components of the process which lead the path into desired regions with positive probability. This is inspired by Chapter~5 of \cite[p.148]{Sato1999}. In order to do so we will have to make a case distinction since the suitable components of the process differ from case to case. We will work with the following general decomposition, which then is specified in the case distinctions. Let $(X_t)_{t\geq 0}$ be a Lévy process with characteristic triple $(a,\sigma^2 , \Pi )$ and $X_0 = 0$. Let $\Pi_1$ and $\Pi_2$ be measures on $\R$ such that $\Pi = \Pi_1 + \Pi_2$. Let $\eta \in (0,1)$. We decompose formally
\begin{align}\label{support decomposition}
X_t = Y_t - a' t + \sigma B_t + P_t^\eta + M_t^\eta,
\end{align}
where
\begin{align*}
a' = a + \int_{(-1,1)}x \Pi_1 (\intdiff x) + \int_{(-1,1) \setminus (-\eta , \eta)} x \Pi_2 (dx)
\end{align*}
and $(B_t)_{t\geq 0}$ is a standard Brownian motion and $(Y_t)_{t\geq 0}$, $(P_t^\eta)_{t \geq 0}$ and $(M_t^\eta)_{t\geq 0}$ are Lévy processes such that
\begin{align*}
-\log \left( \EW{e^{i \theta Y_1}} \right) &= \int_{\R } ( 1 - e^{i\theta x}  )\Pi_1 (\intdiff x)
\intertext{and}
-\log \left( \EW{e^{i \theta P_1^\eta}} \right) &= \int_{\R \setminus (-\eta , \eta)} ( 1 - e^{i\theta x}  )\Pi_2 (\intdiff x)
\intertext{and}
-\log \left( \EW{e^{i \theta M_1^\eta}} \right) &= \int_{(-\eta , \eta ) } ( 1 - e^{i\theta x} + i\theta x ) \Pi_2 (\intdiff x).
\end{align*}
Note that $P_t^\eta$ is a compound Poisson process and $M_t^\eta$ is a zero-mean square-integrable martingale with $\EW{M_t^\eta)^2} = t \int_{(-\eta , \eta)} x^2 \Pi (dx)$. This means that by Doob's inequality for every $t,C>0$ we can choose $\eta >0$ such that
\begin{align}\label{doob inequality}
\PR{\sup_{s\leq t} \abs{M_s^\eta } \geq C} \leq \frac{t}{C} \int_{(-\eta , \eta)} x^2 \Pi_2 (dx)<1.
\end{align}

For treating (3) in Proposition~\ref{levy process supportbedingung} we will use the following auxiliary lemma.

\begin{lemma}\label{Levy process support constant boundary lemma}
Let $(X_t)_{t\geq 0}$ be a Markov process and $\mu$ a probability measure on $\R$. Assume that for any $t>0$ and for any $K> x>0$ we have that
\begin{align*}
 \supp ( \PRx{x}{X_t \in \cdot,\tau_K > t} ) = (-\infty , K].
\end{align*}
Then we have for any boundary function $b:[0,\infty] \to [-\infty , \infty]$ and $t>0$ with $\PRx{\mu}{\tau_b > t} > 0$ that
\begin{align*}
\suppex (\PRbx{\mu}{X_t \in \cdot}{\tau_b > t}) = (-\infty , b(t)].
\end{align*}
\end{lemma}
\begin{proof}
Let $b$ and $t>0$ as in the statement. We abbreviate $\Pr = \Pr_\mu$. It holds $b(t) > -\infty$ since $\PR{\tau_b > t}> 0$. Define for $ 0 < \delta < t$ 
\begin{align*}
 K_\delta \coloneqq \inf_{s\in [t-\delta , t]} b(s).
\end{align*}
Furthermore, since $b$ is lower semicontinuous and $\PR{\tau_b > t} >0$, we have for $0< r < t-\delta$ that
\begin{align*}
K_1 \coloneqq \inf_{s \in [r,t-\delta ]} b(s) = \min_{s \in [r,t-\delta ]} b(s) > -\infty.
\end{align*}
Let $s\in [r, t-\delta ]$ such that $b(s) = K_1$. Define $\mu_s \coloneqq \PR{X_s \in \cdot ,\tau_b > s}$. Note that
\begin{align*}
\emptyset \neq \supp (\mu_s ) \subseteq (-\infty , b(s)] = (-\infty , K_1]
\end{align*}
but $\mu_s (\{K_1\})=0$. We can write
\begin{align*}
&\mu_{t-\delta}  \coloneqq  \PRx{\mu_s}{ X_{t-\delta - s} \in \cdot\; , \tau_{K_1} > t-\delta -s}\\
&= \int_{(-\infty , K_1)} \PRx{x}{ X_{t-\delta - s} \in \cdot\; , \tau_{K_1} > t-\delta -s} \mu_s ( \intdiff x).
\end{align*}
The assumption of the statement ensures that
\begin{align*}
(-\infty,K_1) \subseteq \supp (\PRx{x}{ X_{t-\delta - s} \in \cdot\; , \tau_{K_1} > t-\delta -s})
\end{align*}
for every $x< K_1$, which implies
\begin{align*}
  \suppex (\mu_{t-\delta}) = (-\infty , K_1]
\end{align*}
with $\mu_{t-\delta} (\{K_1\}) =0$. Now let $z\in (-\infty , K_\delta )$ and  $\varepsilon \in (0, K_\delta - z)$. Note that due to the assumption of the statement we have that
\begin{align*}
(-\infty,K_\delta) \subseteq \supp (\PRx{x}{X_\delta \in \cdot, \tau_{K_\delta} > \delta } )
\end{align*}
for every $x< K_\delta$. Thus, using the Markov property, we have
\begin{align*}
&\PR{X_t \in (z-\varepsilon , z + \varepsilon) , \tau_b >t} \\
&\geq \PR{ X_t \in (z-\varepsilon , z + \varepsilon) , \tau_b > s , \tau_{K_1} \notin [s, t-\delta ] , \tau_{K_\delta} \notin [t-\delta , t]}\\
&= \int_{(-\infty , K_1)} \PRx{x}{X_\delta \in (z-\varepsilon , z + \varepsilon), \tau_{K_\delta} > \delta }  \mu_{t-\delta} ( \intdiff x) >0.
\end{align*}
This means that 
\begin{align*}
(-\infty,K_\delta) \subseteq \supp (\PR{X_t \in \cdot\; , \tau_b >t}).
\end{align*}
But since $K_\delta \to b(t)$ as $\delta \to 0$, we have that
\begin{align*}
\suppex (\PR{X_t \in \cdot\; , \tau_b >t}) = (-\infty , b(t)].
\end{align*}
This completes the proof.
\end{proof}

Let us establish conditions which imply the condition of the auxiliary lemma.

\begin{lemma}\label{Levy process support constant boundary}
Let $(X_t)_{t\geq 0}$ be a Lévy process with a characteristic triple $(a,\sigma^2 , \Pi )$ and $X_0=0$. If one of the following holds,
\begin{itemize}
\item[(i)] $(X_t)_{t\geq 0}$ has unbounded variation, i.e. $\sigma \neq 0$ or $\int_\R (1\wedge \abs{x}) \Pi (\intdiff x) = \infty$,
\item[(ii)] $0\in \supp (\Pi)$ and $\Pi ((-\infty , 0))>0$ and $\Pi ( (0 , \infty ))>0$,
\end{itemize}
then
\begin{align*}
\supp (\PRb{X_t \in \cdot\;}{\tau_K>t} ) = (-\infty,K].
\end{align*}
\end{lemma}

We will prove this lemma by using components which have the following form.

\begin{lemma}\label{Levy process support constant boundary hilfslemma}
Let $(X_t)_{t\geq 0}$ be a Lévy process with a characteristic triple $(a,\sigma^2 , \Pi )$ and $X_0=0$. If $\sigma = 0$, $\Pi (\R) < \infty$ and for $\gamma \coloneqq -\left(a + \int_{(-1,1)} x \Pi (\intdiff x)\right)$ one of the following conditions is fulfilled,
\begin{enumerate}
\item[(a)] $\gamma \leq 0$ and $0 \in \supp (\Pi(\cdot \cap (0,\infty)))$,
\item[(b)] $\gamma \geq 0$ and $0 \in \supp (\Pi(\cdot \cap (-\infty,0)))$,
\item[(c)] $0\in \supp (\Pi)$ and $\Pi ((-\infty , 0 ))>0$ and $ \Pi ((0 , \infty ))>0$,
\end{enumerate}
then
\begin{align*}
\supp (\PRb{X_t \in \cdot\;}{\tau_K>t} ) = 
\begin{cases} 
(\gamma t, K] &: \text{ (a)},\\
(-\infty , \min( K , \gamma t)] &: \text{ (b)},\\
(-\infty , K] &: \text{ (c)}.
\end{cases}
\end{align*}
\end{lemma}

\begin{proof}
For $\Pi_1 \coloneqq \Pi$ and $\gamma = -a'$ the decomposition of \eqref{support decomposition} reduces to
\begin{align*}
X_t = Y_t + \gamma t
\end{align*}
For $c_1 < c_2$ and $\kappa_1 < \kappa_2$ let us define
\begin{align*}
Y_t^c \coloneqq \int_0^t (X_s - X_{s-})\intdiff N_s^c , \qquad Y_t^\kappa \coloneqq \int_0^t (X_s - X_{s-})\intdiff N_s^\kappa,
\end{align*}
into independent processes, where
\begin{align*}
N_t^c \coloneqq \sum_{s\leq t} \If_{(c_1 , c_2)}(X_s - X_{s-}), 
\qquad
N_t^\kappa \coloneqq \sum_{s\leq t} \If_{(\kappa_1 , \kappa_2)}(X_s - X_{s-})
\end{align*}
are Poisson processes with intensities $\Pi ((c_1, c_2))$ and $\Pi ((\kappa_1 , \kappa_2))$, respectively.

\textbf{Assume condition (a).} Let $x \in (\gamma t,  K )$ and let $\varepsilon > 0$ such that
\begin{align*}
\varepsilon < \min\{ x -\gamma t ,  K - x \}.
\end{align*}
Since $0 \in \supp (\Pi (\cdot \;\cap (0,\infty ) ))$ there is $\kappa \in (0, \varepsilon /2 ) \cap \supp (\Pi)$. Let $\delta \in (0,t)$ so that
\begin{align*}
-\gamma \delta < \frac{\varepsilon}{2}.
\end{align*}
Since $\kappa < \varepsilon /2$ and $\gamma (t-\delta ) < x -\varepsilon /2$ there is $n_\kappa \in \N$ such that
\begin{align*}
\gamma (t-\delta ) + n_\kappa\cdot \kappa \in \left( x- \frac{\varepsilon}{2} , x+ \frac{\varepsilon}{2} \right)
\end{align*}
There are $0 < \kappa_1 < \kappa < \kappa_2 < \varepsilon / 2$ such that
\begin{align*}
\gamma (t-\delta ) + n_\kappa \cdot(\kappa_1 , \kappa_2) \subseteq \left( x- \frac{\varepsilon}{2} , x+ \frac{\varepsilon}{2} \right)
\end{align*}
Since $\kappa \in \supp (\Pi)$ we have that $\Pi ( (\kappa_1 ,\kappa_2 )) > 0$. With the decomposition
\begin{align*}
Y_t = Y_t^0 + Y_t^\kappa
\end{align*}
we observe that
\begin{align*}
&\left\{Y_{s} = 0\; \forall s \leq t-\delta \right\} \cap \left\{Y_s^0 = 0 \; \forall s \leq t \right\} \cap \left\{ N_t^\kappa = n_\kappa \right\}\\
&\subseteq \left\{ Y_t +\gamma t  \in \left( x-\varepsilon , x+\varepsilon \right) \right \} \cap \left\{ \sup_{s\leq t} (Y_s +\gamma s) <  K \right\}.
\end{align*}
By independence, the Markov property, the fact that the intensity of $(Y_s)_{s\geq 0}$ and $(Y_s^0)_{s\geq 0}$ is finite and that $\Pi ((\kappa_1, \kappa_2))>0$, we have that
\begin{align*}
&\PR{Y_{s} = 0\; \forall s \leq t-\delta , Y_s^0 = 0 \; \forall s \leq t , N_t^\kappa = n_\kappa}\\
&= \PR{Y_{s} = 0\; \forall s \leq t-\delta} \PR{Y_s^0 = 0 \; \forall s \leq \delta} \PR{N_\delta^\kappa = n_\kappa} >0 .
\end{align*}
This means that in the situation of (a) we have that
\begin{align*}
\supp ( \PR{X_t \in \cdot \; , \tau_K > t} ) = [\gamma t , K].
\end{align*}

\textbf{Assume condition (b).} Let $x\in (-\infty , \min(K, \gamma t))$ and let $\varepsilon <  \min(K, \gamma t)-x$. Let $\delta \in (0,t)$ such that
\begin{align*}
\gamma \delta < K. 
\end{align*}
Since $0\in \supp ( \Pi ( \cdot \cap (-\infty , 0) ))$ there is $\kappa \in (-\varepsilon , 0) \cap \supp (\Pi)$. Since $\abs{\kappa } < \varepsilon$ and $x+ \varepsilon < \gamma t$ there is $n_\kappa \in \N$ such that
\begin{align*}
\gamma t + n_\kappa \cdot \kappa \in \left( x- \varepsilon , x+ \varepsilon \right)  .
\end{align*}
There are $-\varepsilon < \kappa_1 < \kappa < \kappa_2 <0$ such that
\begin{align*}
\gamma t + n_\kappa \cdot (\kappa_1 , \kappa_2) \subseteq \left( x- \varepsilon , x+ \varepsilon \right)  .
\end{align*}
Since $\kappa \in \supp (\Pi)$ we have that $\Pi ( (\kappa_1 ,\kappa_2 )) > 0$. With the decomposition
\begin{align*}
Y_t = Y_t^0 + Y_t^\kappa
\end{align*}
we observe that
\begin{align*}
&\left\{Y_s^0 = 0 \; \forall s \leq t \right\} \cap \left\{ N_\delta^\kappa = n_\kappa \right\} \cap \left\{ N_s^\kappa = n_\kappa \;\forall \delta \leq s \leq t \right\}\\
&\subseteq\left\{ Y_t +\gamma t \in \left( x- \varepsilon , x+ \varepsilon \right), \sup_{s\leq t} (Y_s +\gamma s ) < K \right\}.
\end{align*}
By independence, the Markov property, the fact that the intensity of $(Y_s)_{s\geq 0}$ and $(Y_s^0)_{s\geq 0}$ is finite and that $\Pi ((\kappa_1, \kappa_2))>0$, we have that
\begin{align*}
&\PR{Y_s^0 = 0 \; \forall s \leq t ,  N_\delta^\kappa = n_\kappa , N_s^\kappa = n_\kappa \;\forall \delta \leq s \leq t }\\
&= \PR{Y_s^0 = 0 \; \forall s \leq t}\PR{ N_\delta^\kappa = n_\kappa}\PR{N_s^\kappa = 0 \;\forall 0 \leq s \leq t -\delta} >0.
\end{align*}
This means that in the situation of (b) we have that
\begin{align*}
\supp ( \PR{X_t \in \cdot \; , \tau_K > t} ) = (-\infty , \min (K , \gamma t)].
\end{align*}

\textbf{Assume the case (c).} We  have that
\begin{align*}
\text{(c.1) } 0 \in \supp ( \Pi (\cdot \cap (0,\infty))) \quad \text{ or } \quad \text{(c.2) } 0\in \supp (\Pi (\cdot \cap (-\infty , 0))).
\end{align*}
\textbf{Suppose that (c.1) holds:} Let $x\in (-\infty , K)$ and $\varepsilon >0$ such that
\begin{align*}
\varepsilon < K - x.
\end{align*}
Let $\delta_1 ,\delta_2 \in (0,t)$ such that
\begin{align*}
\gamma \delta_1 < K \quad \text{and} \quad \abs{\gamma \delta_2} < \frac{\varepsilon}{2}.
\end{align*}
By condition (c) there is $c\in (-\infty , 0) \cap \supp (\Pi)$. Thus there is $m_c\in \N$ such that
\begin{align*}
m_c \cdot c +\gamma (t-\delta_2) < x-\frac{\varepsilon}{2}.
\end{align*}
Since $0 \in \supp (\Pi (\cdot \;\cap (0,\infty ) ))$ there is $\kappa \in (0, \varepsilon ) \cap \supp (\Pi)$. Since $\kappa < \varepsilon $ and $m_c \cdot c + \gamma (t-\delta_2) < x-\varepsilon /2$ there is $n_\kappa \in \N$ such that
\begin{align*}
m_c \cdot c +\gamma (t-\delta_2) + n_\kappa\cdot \kappa \in \left( x- \frac{\varepsilon}{2}, x+ \frac{\varepsilon}{2} \right)
\end{align*}
There are $c_1 < c < c_2 < 0$ and $0 < \kappa_1 < \kappa < \kappa_2 < \varepsilon $ such that
\begin{align*}
m_c \cdot (c_1,c_2) + \gamma (t-\delta_2) + n_\kappa \cdot(\kappa_1 , \kappa_2) \subseteq \left( x- \frac{\varepsilon}{2} , x+ \frac{\varepsilon}{2} \right)
\end{align*}
Since $c,\kappa \in \supp (\Pi)$ we have that $\Pi ( (c_1 ,c_2 )) > 0$ and $\Pi ( (\kappa_1 ,\kappa_2 )) > 0$. Now observe that
\begin{align*}
&\left\{ Y_s^0 = 0 \; \forall s \leq t \right\} \cap\{N_s^c = m_c \;\forall \delta_1 \leq s \leq t\}\cap \{N_s^\kappa = 0 \;\forall s \leq t -\delta_2\}\cap \{N_t^\kappa = n_\kappa\}\\
&\subseteq \left\{ Y_t +\gamma t \in (x-\varepsilon , x+\varepsilon) , \sup_{s\leq t} (Y_s + \gamma s) < K\right\}.
\end{align*}
By independence, the Markov property and the fact that the intensities of $(Y_s^0)_{s\geq 0}$, $(Y_s^\eta)_{s\geq 0}$ and $(Y_s^c)_{s\geq 0}$ are finite and for $(Y_s^\eta)_{s\geq 0}$ and $(Y_s^c)_{s\geq 0}$ even positive, we have that
\begin{align*}
&\PR{ Y_s^0 = 0 \; \forall s \leq t , N_s^c = m_c \;\forall \delta_1 \leq s \leq t, N_s^\kappa = 0 \;\forall s\leq t -\delta_2, N_t^\kappa = n_\kappa}\\
&= \PR{ Y_s^0 = 0 \; \forall s \leq t} \PR{N_{\delta_1}^c = m_c}\PR{N_s^c = 0 \;\forall  s \leq t-\delta_1}\\
&\phantom{= \PR{ Y_s^0 = 0 \; \forall s \leq t}}\cdot \PR{N_s^\kappa = 0 \;\forall s \leq t -\delta_2 }\PR{N_{\delta_2}^\kappa = n_\kappa} >0
\end{align*}
This finishes the proof of the lemma for the case (c.1).

\textbf{Suppose that (c.2) holds:} Let $x\in (-\infty , K)$ and $\varepsilon >0$ such that
\begin{align*}
\varepsilon < K - x.
\end{align*}
Let $\delta_1 ,\delta_2 \in (0,t)$ such that
\begin{align*}
\gamma \delta_1 < K \quad \text{and}\quad \abs{\gamma \delta_2} < \frac{\varepsilon}{2}.
\end{align*}
By condition (c) there is $c\in (0 , \infty) \cap \supp (\Pi)$. Thus there is $m_c\in \N$ such that
\begin{align*}
m_c \cdot c +\gamma (t-\delta_2) > x+\frac{\varepsilon}{2}.
\end{align*}
Since $0 \in \supp (\Pi (\cdot \;\cap (-\infty,0 ) ))$ there is $\kappa \in (-\varepsilon ,0) \cap \supp (\Pi)$. Since $\abs{\kappa} < \varepsilon $ and $m_c \cdot c + \gamma (t-\delta_2) > x+\varepsilon /2$ there is $n_\kappa \in \N$ such that
\begin{align*}
m_c \cdot c + \gamma (t-\delta_2) + n_\kappa\cdot \kappa \in \left( x- \frac{\varepsilon}{2}, x+ \frac{\varepsilon}{2} \right)
\end{align*}
There are $0<c_1 < c < c_2 $ and $- \varepsilon< \kappa_1 < \kappa < \kappa_2 < 0 $ such that
\begin{align*}
m_c \cdot (c_1,c_2) + \gamma (t-\delta_2) + n_\kappa \cdot(\kappa_1 , \kappa_2) \subseteq \left( x- \frac{\varepsilon}{2} , x+ \frac{\varepsilon}{2} \right)
\end{align*}
Since $c,\kappa \in \supp (\Pi)$ we have $\Pi ( (c_1 ,c_2 )) > 0$ and $\Pi ( (\kappa_1 ,\kappa_2 )) > 0$. By decomposing
\begin{align*}
Y_t = Y_t^0 + Y_t^c + Y_t^\kappa
\end{align*}
we observe that
\begin{align*}
&\left\{ Y_s^0 = 0 \; \forall s \leq t \right\} \cap\{N_s^\kappa = n_\kappa \;\forall \delta_1 \leq s \leq t\}\cap \{N_s^c = 0 \;\forall s \leq t -\delta_2\}\cap \{N_t^c = m_c\}\\
&\subseteq \left\{ Y_t +\gamma t \in (x-\varepsilon , x+\varepsilon) , \sup_{s\leq t} (Y_s +\gamma s) < K\right\}.
\end{align*}
By independence, the Markov property and the fact that the intensities of $(Y_s^0)_{s\geq 0}$, $(Y_s^\kappa)_{s\geq 0}$ and $(Y_s^c)_{s\geq 0}$ are finite and for $(Y_s^\kappa)_{s\geq 0}$ and $(Y_s^c)_{s\geq 0}$ even positive, we have
\begin{align*}
&\PR{ Y_s^0 = 0 \; \forall s \leq t , N_s^\kappa = n_\kappa \;\forall \delta_1 \leq s \leq t, N_s^c = 0 \;\forall s\leq t -\delta_2, N_t^c = m_c}\\
&= \PR{ Y_s^0 = 0 \; \forall s \leq t} \PR{N_{\delta_1}^\kappa = n_\kappa}\PR{N_s^\kappa = 0 \;\forall  s \leq t-\delta_1}\\
&\phantom{= \PR{ Y_s^0 = 0 \; \forall s \leq t}}\cdot \PR{N_s^c = 0 \;\forall s \leq t -\delta_2 }\PR{N_{\delta_2}^c = m_c} >0
\end{align*}
This finishes the proof of (c.2) and thus for the situation of (c).
\end{proof}

\begin{proof}[Proof of Lemma~\ref{Levy process support constant boundary}]
\textbf{Let us consider the case $\sigma >0$.} For $\Pi_2 \coloneqq \Pi$ and $\eta \in (0,1)$ the decomposition of \eqref{support decomposition} reads
\begin{align*}
X_t = - a_\eta t + \sigma B_t + P_t^\eta + M_t^\eta
\end{align*}
with $a_\eta = a + \int_{(-1,1) \setminus (-\eta , \eta)} x \Pi (dx)$.

Now let $x \in (-\infty , K)$. Let $\varepsilon >0$ be such that
\begin{align*}
\varepsilon < \min \left\{ K-x , K  \right\}.
\end{align*}
Due to \eqref{doob inequality} we can choose $\eta > 0$ such that
\begin{align*}
\PR{\sup_{s\leq t} \Abs{ M_s^\eta } \geq \frac \varepsilon 2 } < 1.
\end{align*}
Further we have that $\PR{P_{s}^\eta = 0\; \forall s \leq t} > 0$. Let $f: [0,t] \to \R$ be defined by $f(s) \coloneqq \frac{s}{t}x$. From the theory of Brownian motion we know that
\begin{align*}
\PR{ \Abs{ \sigma B_s - a_\eta s - f(s)} < \frac \varepsilon 2 \;\forall s \in [0,t]}>0,
\end{align*}
for example see Theorem 38 in \cite{Freedman1983}. Note that
\begin{align*}
&\left\{\Abs{ \sigma B_s - a_\eta s - f(s)} < \frac \varepsilon 2 \;\forall s \in [0,t] \right\}\\
&\subseteq \left\{ \Abs{\sigma B_t - a_\eta t - x} < \frac \varepsilon 2 \right\} \cap \left\{ \sup_{s\leq t}\Abs{ \sigma B_s - a_\eta s} < \frac \varepsilon 2 + \max \left\{ 0 , x \right\} \right\}.
\end{align*}
This means
\begin{align*}
&\left\{\Abs{ \sigma B_s - a_\eta s - f(s)} < \frac \varepsilon 2 \;\forall s \in [0,t] \right\} \cap \left\{ P_s^\eta = 0 \; \forall s \leq t \right\} \cap \left\{ \sup_{s\leq t} \Abs{ M_s^\eta } < \frac \varepsilon 2 \right\}\\
&\subseteq \left\{ \Abs{X_t -x } < \varepsilon \right\} \cap \left\{ \sup_{s\leq t} X_s < K \right\}.
\end{align*}
This yields
\begin{align*}
0 &< \PR{\Abs{ \sigma B_s - a_\eta s - f(s)} < \frac \varepsilon 2 \;\forall s \leq t}\PR{P_s^\eta = 0 \; \forall s \leq t} \PR{ \sup_{s\leq t} \Abs{ M_s^\eta} < \frac{\varepsilon}{2} } \\
&\leq \PR{X_t \in (x-\varepsilon, x+\varepsilon ) ,  \sup_{s\leq t} X_s < K } = \PR{X_t \in (x-\varepsilon , x+\varepsilon ), \tau_K > t}.
\end{align*}
Thus we have $\supp(\PR{X_t \in \cdot\;, \tau_K > t}) = (-\infty , K]$.

\textbf{Now assume that $\sigma^2=0$}. Define $\Pi_1 (\intdiff x) \coloneqq x^2 \If_{(-1,1)}(x)\Pi(\intdiff x)$ and $\Pi_2 \coloneqq \Pi - \Pi_1$. Observe that it holds $\Pi_1 (\R)<\infty$. For $\eta \in (0,1)$ the decomposition of \eqref{support decomposition} reads
\begin{align*}
X_t = Y_t - a_\eta t + P_t^\eta + M_t^\eta
\end{align*}
with $a_1 \coloneqq a + \int_\R x \Pi_1 (\intdiff x)$ and $a_\eta = a_1 + \int_{(-1,1)\setminus (-\eta,\eta)} x \Pi_2 (dx)$.

In the following we distinguish the following three cases.
\begin{itemize}
\item[(I)] $\Pi ((-\infty , 0))=0$ and $\int_\R (1 \wedge \abs{x} ) \Pi (\intdiff x) = \infty$,
\item[(II)] $\Pi ( (0 , \infty ))=0$ and $ \int_\R (1 \wedge \abs{x} ) \Pi (\intdiff x) = \infty$,
\item[(III)] $0\in \supp (\Pi)$ and $\Pi ( ( -\infty , 0 ))>0$ and $\Pi ( (0 , \infty ))>0$.
\end{itemize}
Note that these cases are exhausting for (i) and (ii) if $\sigma^2 =0$.

Let $t>0$ and $K>0$. Let $x\in (-\infty , K)$ and $\varepsilon > 0$ such that
\begin{align*}
\varepsilon < \min \{K-x, K\}.
\end{align*}

We first claim that for all three cases there exists $\eta>0$ such that
\begin{align*}
\PR{\sup_{s\leq t } \abs{M_s^\eta } \geq  \frac{\varepsilon}{2}} <1
\end{align*}
and
\begin{align*}
\PR{ \Abs{Y_t - ta_\eta -x} < \frac{\varepsilon}{2}  , \sup_{s\leq t} (Y_s -sa_\eta) < K- \frac{\varepsilon}{2}}>0.
\end{align*}

Let us assume for the moment that the claim is true. We will finish the proof of the theorem from here and prove the claim further below. Recall that $X_t = Y_t - a_\eta t + P_t^\eta + M_t^\eta$ , thus
\begin{align*}
&\left\{ \Abs{Y_t - ta_\eta -x} < \frac{\varepsilon}{2} , \sup_{s\leq t} (Y_s -sa_\eta) < K- \frac{\varepsilon}{2} \right\} \cap \{ P_s^\eta = 0 \; \forall s \leq t \} \cap \left\{ \sup_{s\leq t } \abs{M_s^\eta } < \frac{\varepsilon}{2} \right\}\\
&\subseteq \left\{ \Abs{ X_t -x} < \varepsilon , \sup_{s\leq t} X_s < K \right\} = \{ X_t \in (x-\varepsilon , x+ \varepsilon ) , \tau_K > t\}.
\end{align*}
Note that $\PR{P_s^\eta = 0 \; \forall s \leq t}>0$ and $\PR{\sup_{s\leq t } \abs{M_s^\eta } < \frac{\varepsilon}{2}} >0$ and by independence of the events we obtain
\begin{align*}
\PR{ X_t \in (x-\varepsilon , x+ \varepsilon ) , \tau_K > t} > 0,
\end{align*}
which implies the statement of the theorem, since $x\in (-\infty , K)$ was arbitrary.

Let us now prove the claim for every case separately.\\
\textbf{Let us assume (I).} This implies that for all $\eta >0$ we have $\Pi ((0,\eta ) ) = \infty$, hence
\begin{align*}
0 \in \supp (\Pi(\cdot \cap (0,\infty))).
\end{align*}
Moreover, this implies that for $\eta \in (0, 1)$ we have $a_\eta = a_1 + \int_{[\eta ,1)} x \Pi_2 (dx)$ and
\begin{align*}
\int_{(0,1)} x \Pi_2 (dx)  = \infty,
\end{align*}
which implies $\lim_{\eta \to 0 } a_\eta = \infty$. With \eqref{doob inequality} in mind we can choose $\eta >0$ such that
\begin{align*}
a_\eta t\geq \max\{- (x-\varepsilon /2) ,0\} \quad \text{ and } \quad \PR{\sup_{s\leq t } \abs{M_s^\eta } \geq \frac{\varepsilon}{2}} < 1.
\end{align*}
Observe that, since $0\in \supp (\Pi_1 (\cdot \cap (0,\infty)))$ and $a_\eta \geq 0$, the condition (a) of Lemma~\ref{Levy process support constant boundary hilfslemma} is fulfilled for the process $\tilde{X}_t = Y_t -a_\eta t$, thus we have that
\begin{align*}
\supp \left( \PR{Y_t - a_\eta t \in \cdot \; , \sup_{s\leq t} (Y_s -a_\eta s) < K- \frac{\varepsilon}{2}} \right) = \left( -a_\eta t , K- \frac{\varepsilon}{2} \right],
\end{align*}
which implies the assertion of the claim. This finishes the proof for the case (I).

\textbf{Let us assume (II).} Then for all $\eta >0$ we have $\Pi ((-\eta ,0 ) ) = \infty$, hence
\begin{align*}
0 \in \supp (\Pi(\cdot \cap (-\infty,0))).
\end{align*}
Moreover, this implies that for $\eta \in (0, 1)$ we have $a_\eta = a_1 + \int_{(1,-\eta]} x \Pi_2 (dx)$ and
\begin{align*}
\int_{(-1,0)} x \Pi_2 (dx)  = -\infty,
\end{align*}
which implies $\lim_{\eta \to 0 } a_\eta = -\infty$. With \eqref{doob inequality} in mind we can choose $\eta >0$ such that
\begin{align*}
- a_\eta t\geq K - \frac{\varepsilon}{2} \quad \text{ and } \quad \PR{\sup_{s\leq t } \abs{M_s^\eta } \geq \frac{\varepsilon}{2}} < 1.
\end{align*}
Observe that now condition (b) of Lemma~\ref{Levy process support constant boundary hilfslemma} is fulfilled for the process $\tilde X_t = Y_t - a_\eta t$, thus we have that
\begin{align*}
\supp \left( \PR{Y_t - a_\eta t \in \cdot \; , \sup_{s\leq t} (Y_s -a_\eta s) < K- \frac{\varepsilon}{2}} \right) = \left( -\infty , K- \frac{\varepsilon}{2} \right],
\end{align*}
which implies the assertion of the claim. This finishes the proof for the case (II).

\textbf{Let us assume (III).} Due to \eqref{doob inequality} we can choose $\eta > 0$ such that
\begin{align*}
\PR{\sup_{s\leq t } \abs{M_s^\eta } \geq \frac{\varepsilon}{2}} < 1.
\end{align*}
The process $\tilde X_t = Y_t - a_\eta t$ inherits the properties of (ii) and fulfills the conditions of (c) of Lemma~\ref{Levy process support constant boundary hilfslemma}. Thus we have that
\begin{align*}
\supp \left( \PR{Y_t - a_\eta t \in \cdot \; , \sup_{s\leq t} (Y_s -a_\eta s) < K- \frac{\varepsilon}{2}} \right) = \left( -\infty , K- \frac{\varepsilon}{2} \right].
\end{align*}
This implies the assertion of the claim. This finishes the proof for the case (III).
\end{proof}

\begin{proof}[Proof of Proposition~\ref{levy process supportbedingung}]
\textbf{Regarding (III)} note that by translation of the starting point Lemma~\ref{Levy process support constant boundary} yields that for any $t>0$ and for any $K>x>0$ we have that
\begin{align*}
 \supp ( \PRx{x}{X_t \in \cdot,\tau_K > t} ) = (-\infty , K].
\end{align*}
Thus we obtain by Lemma~\ref{Levy process support constant boundary lemma} that for any boundary function $b:[0,\infty] \to [-\infty , \infty]$ and $t>0$ with $\PR{\tau_b > t} > 0$ it holds
\begin{align*}
\suppex (\PRb{X_t \in \cdot}{\tau_b > t}) = (-\infty , b(t)].
\end{align*}
This finishes the proof for (3).

\textbf{Let us prove (I) and (II)}: For $-1 < \kappa_1 < \kappa_2<1$ define $\Pi_1 (\intdiff x) \coloneqq \abs{x}\If_{(\kappa_1,\kappa_2)} (x) \Pi (\intdiff x)$ and $\Pi_2 \coloneqq \Pi -\Pi_1$. For $\eta \in (0,1)$ the decomposition of \eqref{support decomposition} reads
\begin{align*}
X_t = Y_t + P^\eta_t + S^\eta_t,
\end{align*}
where $S_t^\eta \coloneqq M_t^\eta - a_\eta t$ and since there is no drift, i.e. $a_\eta = \int_{(-\eta , \eta)} x \Pi_2 (\intdiff x)$,
\begin{align*}
-\log \left( \EW{e^{i \theta S_1^\eta}} \right) &= \int_{(-\eta , \eta ) } ( 1 - e^{i\theta x}) \Pi_2 (\intdiff x).
\end{align*}
Let us denote
\begin{align*}
\tau^\eta_b \coloneqq \inf\{ t > 0 : S_t^\eta \geq b(t) \}.
\end{align*}
Since there is no drift, in both of the cases (1) or (2) the process $( \abs{S_t^\eta })_{t\geq 0}$ is a subordinator with $\EW{ S_t^\eta} = t \int_{(-\eta , \eta)} \abs{x} \Pi_2 (dx)$. This means that by Markov's inequality for every $C>0$ we can choose $\eta >0$ such that
\begin{align}\label{Markov inequality}
\PR{\sup_{s\leq t} \abs{S_s^\eta } \geq C} = \PR{ \abs{S_t^\eta } \geq C}\leq \frac{t}{C} \int_{(-\eta , \eta)} \abs{x} \Pi_2 (dx)<1.
\end{align}
For $\kappa_1 < \kappa_2$ we write
\begin{align*}
Y_t = \int_0^t (Y_s-Y_{s-}) \intdiff N_s^\kappa \quad \text{with}\quad 
N_t^\kappa &\coloneqq \sum_{s\leq t} \If_{\R \setminus \{0\}}(Y_s - Y_{s-}),
\end{align*}
and $N_t^\kappa$ is a Poisson process with rate $\Pi_1 ((\kappa_1 ,\kappa_2))\geq 0$.

Now we treat (1) and (2) separately.\\
\textbf{Assume the conditions of (I).} Assume that $\PR{\tau_b>t} >0$ for $t>0$. We have $\supp (\Pi) \subseteq [0,\infty)$. Thus we have $X_r \geq N_r^\eta $ almost surely. This implies
\begin{align*}
0 < \PR{\tau_b>t} \leq \PR{\tau_b^\eta > t}.
\end{align*}
Let $\delta \in (0,t)$ and define
\begin{align*}
K_\delta \coloneqq \inf_{s\in [t-\delta , t] } b(s) > 0.
\end{align*}
Let $x\in (0,K_\delta)$ and $0< \varepsilon < \min (x , K_\delta -x)$. By \eqref{Markov inequality} choose $\eta \in (0,1)$ such that
\begin{align*}
\PR{\sup_{s\leq t} \abs{S_s^\eta } \geq \frac{2\varepsilon}{3}} < \PR{\tau_b>t}.
\end{align*}
Since $ \tau_b^\eta \geq \tau_b$ almost surely, we have that
\begin{align*}
&\PR{\sup_{s\leq t} \abs{S_s^\eta } < \frac{2\varepsilon}{3} , \tau_b^\eta > t}\geq \PR{\sup_{s\leq t} \abs{S_s^\eta } < \frac{2\varepsilon}{3} , \tau_b > t}\\
&\geq \PR{\tau_b > t} -\PR{\sup_{s\leq t} \abs{S_s^\eta } \geq \frac{2\varepsilon}{3}} >0.
\end{align*}
Further, since $0\in \supp(\Pi) \subseteq [0,\infty)$ we have that there is $\kappa \in \supp(\Pi) \cap (0,\varepsilon/3)$. Since $\kappa < \varepsilon /3$ there is $n_\kappa \in \N$ and $0< \kappa_1 < \kappa < \kappa_2 < \varepsilon /3$ such that
\begin{align*}
(\kappa_1 ,\kappa_2) \cdot n_\kappa \subseteq \left( x -\frac{\varepsilon}{3}, x+\frac{\varepsilon}{3}\right). 
\end{align*}
Now observe, since $(Y_s)_{s\geq 0}$ only has jumps of size contained in $(\kappa_1, \kappa_2)$, that
\begin{align*}
&\left\{ \sup_{s\leq t} \abs{S_s^\eta } < \frac{2\varepsilon}{3} , \tau_b^\eta > t \right\} \cap \left\{ P_s^\eta =0 \;\forall s\leq t \right\} \cap \left\{ N_s^\kappa = 0 \; \forall s\leq t-\delta \right\}  \cap \left\{ N_t^\kappa = n_\kappa \right\}\\
&\subseteq \left\{ \sup_{s\leq t} \abs{S_s^\eta} < \frac{2\varepsilon}{3}\right\} \cap \left\{ \tau_b > t-\delta \right\} \cap \left\{ P_s^\eta =0 \;\forall s\leq t \right\}  \\
&\phantom{\subseteq}\quad\cap \left\{ Y_t \in \left( x -\frac{\varepsilon}{3}, x+\frac{\varepsilon}{3}\right) \right\}\\
&\subseteq \left\{ \tau_b > t-\delta \right\}\cap \left\{ \tau_{K_\delta}  \notin [t-\delta ,t] \right\} \cap \left\{ X_t \in (x-\varepsilon , x+\varepsilon) \right\}\\
&\subseteq \left\{ \tau_b > t \right\}\cap \left\{ X_t \in (x-\varepsilon , x+\varepsilon) \right\}.
\end{align*}
By independence and the Markov property the event on the left-hand-side has positive probability and thus we obtain
\begin{align*}
(0,K_\delta) \subseteq \supp (\PR{X_t\in\cdot\; , \tau_b > t}).
\end{align*}
Since $K_\delta \to b(t)$ for $\delta \to 0$, we obtain that
\begin{align*}
\suppex (\PR{X_t\in\cdot\; , \tau_b > t}) = [0,b(t)].
\end{align*}
\textbf{Assume the conditions of (II).}
Let $x\in (-\infty , b(t))$ and $\varepsilon \in (0, b(t)-x)$. By \eqref{Markov inequality} choose $\eta \in (0,1)$ such that
\begin{align*}
\PR{\sup_{s\leq t} \abs{S_s^\eta } < \frac{2\varepsilon}{3}} >0.
\end{align*}
Assume $\PR{\tau_b >0 }>0$. By Blumenthal's law we have that $\tau_b > 0$ almost surely. Now, since $\Pi_1 (\R)<\infty$ and $\Pi_2 ( \R\setminus (-\eta , \eta ))<\infty$, we have
\begin{align*}
0 &< \PR{ Y_s = 0 \;\forall s\leq t, P_s^\eta= 0 \;\forall s\leq t} = \PR{ Y_s = 0 \;\forall s\leq t, P_s^\eta= 0 \;\forall s\leq t, \tau_b > 0}\\
&= \PR{ Y_s = 0 \;\forall s\leq t, P_s^\eta= 0 \;\forall s\leq t, \tau_b^\eta > 0} \\
&= \PR{ Y_s = 0 \;\forall s\leq t, P_s^\eta= 0 \;\forall s\leq t} \PR{ \tau_b^\eta > 0}.
\end{align*}
Consequently we have $ \tau_b^\eta > 0$ almost surely. Hence there is $\delta \in (0,t)$ such that
\begin{align*}
\PR{\sup_{s\leq t} \abs{S_s^\eta } < \frac{2\varepsilon}{3} , \tau_b^\eta > \delta } >0.
\end{align*}
Now since $0\in \supp (\Pi)\cap (-\infty , 0)$ there is $\kappa \in \supp (\Pi) \cap (-\varepsilon , 0 )$. Since $\abs{\kappa}< \varepsilon$ there is $n_\kappa \in \N$ and $-\varepsilon < \kappa_1 < \kappa<\kappa_2 <0$ such that
\begin{align*}
(\kappa_1 , \kappa_2 ) \cdot n_\kappa \subseteq \left( x -\frac{\varepsilon}{3}, x +\frac{\varepsilon}{3}\right).
\end{align*}
Since in the situation of (II) we have $\supp (\Pi) \subseteq (-\infty , 0]$, it follows that $X_r \leq S_r^\eta $. Hence also $\tau_b \geq \tau_b^\eta $. Therefore, with using that $b(t) \leq b(u) $ for all $u\in [0,t]$, we have
\begin{align*}
&\left\{ \sup_{s\leq t} \abs{S_s^\eta } < \frac{2\varepsilon}{3} , \tau_b^\eta > \delta  \right\} \cap \left\{ P_s^\eta =0 \;\forall s\leq t \right\} \cap \left\{ N_\delta^\kappa = n_\kappa \right\}  \cap \left\{ N_\delta^\kappa = n_\kappa\; \forall s \in [\delta , t] \right\}\\
&\subseteq \left\{ \sup_{s\leq t} \abs{S_s^\eta} < \frac{2\varepsilon}{3} \right\}\cap \left\{ \tau_b > \delta  \right\}\cap \left\{ P_s^\eta =0 \;\forall s\leq t \right\}\\
&\phantom{\subseteq}\quad \cap \left\{ Y_s^\kappa < b(t) - \frac{2\varepsilon}{3} \; \forall s\in [\delta , t] \right\}\cap \left\{ Y_t \in \left( x -\frac{\varepsilon}{3}, x+\frac{\varepsilon}{3}\right) \right\}\\
&\subseteq \left\{ \tau_b > \delta  \right\} \cap \left\{ \tau_b \notin [\delta ,t]\right\} \cap \left\{ X_t \in (x-\varepsilon , x+\varepsilon) \right\}\\
&= \left\{ \tau_b > t  \right\} \cap \left\{ X_t \in (x-\varepsilon , x+\varepsilon) \right\}.
\end{align*}
By independence and the Markov property the event on the left-hand-side has positive probability and thus we obtain
\begin{align*}
(-\infty , b(t)) \subseteq \supp (\PR{X_t\in\cdot\; , \tau_b > t}).
\end{align*}
This finishes the proof.
\end{proof}

\begin{proof}[Proof of Theorem~\ref{Levy processes result}]
A Lévy process has right-continuous paths by definition. Furthermore, as a cádlág Feller-process a Lévy process is quasi-left-continuous and a Markov process, see Proposition~7 and Proposition~6 of \cite{bertoin1996}. This gives \ref{existence: continuity assumption} and \ref{uniqueness: markov assumption}. For the order-preservation, note that, by Theorem~1.A.1 of \cite{Shaked2007}, we have $\mu_1 \st\mu_2$ if and only if there exist random variables $Z_i \sim \mu_i$ such that $Z_1 \leq Z_2$. We can choose them independently from $(X_t)_{t\geq 0}$, hence we have that $Z_1 + X_t \leq Z_2 + X_t$ and $\PRx{0}{Z_i + X_t \in \cdot \;} = \PRx{\mu}{X_t \in \cdot\;}$. By Theorem~1.A.1 of \cite{Shaked2007} it follows that $\PRx{\mu_1}{X_t\in\cdot }\st \PRx{\mu_2}{X_t\in \cdot}$. This gives \ref{uniqueness: order preservation assumption}.\\
\textbf{Existence:} By Proposition~\ref{eks first passage time variant} we obtain that (\ref{existence: diffusitivy assumption} $\Rightarrow$ \ref{existence: first-passage time assumption}). Assuming \ref{existence: diffusitivy assumption} therefore implies that the conditions of Theorem~\ref{existence} are fulfilled, and thus a solution for the inverse first-passage time problem exists if \ref{existence: diffusitivy assumption} holds.\\
\textbf{Uniqueness:} It is left to show that in the situation of (a) or (b) we have \ref{uniqueness: support assumption} for the corresponding $I^\xi \subset (0,t^\xi )$. Let $b$ be a boundary function with $\tau_b \overset{\text d} = \xi$. We want to apply Proposition~\ref{levy process supportbedingung}.\\
Assume (a): Let $t\in I^\xi = (0,t^\xi)$. We can exhaust (a) by the case distinction
\begin{enumerate}
\item[(a.i)] $(X_t)_{t\geq 0}$ has unbounded variation,
\item[(a.ii.1)] $0 \in \supp (\Pi)$ and $\Pi ((-\infty , 0))>0$ and $\Pi((0,\infty))>0$,
\item[(a.ii.2')] $0 \in \supp (\Pi)$ and $X_t = Y_t + \gamma t$, where $(Y_t)_{t\geq 0}$ is a subordinator without drift and $\gamma \in \R$.
\end{enumerate}
Note that the case (a.ii.2') can be rephrased as the case that $0 \in \supp (\Pi)$ and $\Pi ((-\infty , 0))=0$ and $(X_t)_{t\geq 0}$ has bounded variation.  Observe, if $x\in \R$, then by a translation according to $x$ and Proposition~\ref{levy process supportbedingung} we have for $t>0$ with $\PRx x{\tau_b > t}>0$ (for (a.ii.2') this implies $x +\gamma t < b(t)$) that
\begin{align*}
\suppex (\PRbx x {X_t \in \cdot \;}{\tau_b > t}) = \begin{cases} (-\infty , b(t)] &: \text{(a.i) or (a.ii.1)} \\ [x+ \gamma t,b(t)] &: \text{(a.ii.2')}. \end{cases}
\end{align*}
Due to $t\in (0,t^\xi)$ we have $ \PR{\tau_b > t}>0 $. Since $\Pr = \int_\R \Pr_x \mu (\intdiff x)$, we obtain
\begin{align*}
\sup \supp (\PRbx \mu {X_t \in \cdot \;}{\tau_b > t}) = b(t).
\end{align*}
Assume (b): Let $t \in I^\xi= \supp (\PR{\xi \in \cdot\;}) \cap (0,t^\xi)$. The case that $(X_t)_{t\geq 0}$ has unbounded variation is already covered by (a). Let us therefore assume that $(X_t)_{t\geq 0}$ has bounded variation. This implies that $X_t = \tilde{X}_t + \gamma t$, where $(-\tilde X_t)_{t\geq 0}$ is a subordinator without drift. Without loss of generality we can assume that $\gamma =0$ by considering the process $(\tilde X_t)_{t\geq 0}$ and the boundary $\tilde b (t) = b(t) - \gamma t$ instead. Hence from now on we assume that $(-X_t)_{t\geq 0}$ is a subordinator without drift and $0\in \supp (\Pi)$.
Suppose there is $0\leq u < t$ such that $b(t) > b(u)$. For $\varepsilon \in (0, b(t) - b(u))$ there exists $\delta \in (0, t-u)$ such that
\begin{align*}
\inf_{s\in [t-\delta , t + \delta]} b(s) \geq b(t) - \varepsilon > b(u).
\end{align*}
But this implies, since $(X_t)_{t\geq 0}$ has non-increasing paths, that
\begin{align*}
0 &< \PR{\xi \in (t-\delta , t+\delta)} = \PR{\tau_b \in (t-\delta , t+\delta)} = \PR{\tau_b \in (t-\delta , t+\delta) , \tau_b > u}  \\
&\leq \PR{\exists s \in (t-\delta , t+\delta) : X_s \geq b(u) , \forall s \geq u : X_s < b(u)} = 0.
\end{align*}
This contradiction shows that $b(t) \leq b(u)$ for all $u\leq t$. For $x\in \R$ with $\PRx x{\tau_b > 0}>0$, by Proposition~\ref{levy process supportbedingung}, we get that
\begin{align*}
\supp (\PRbx x {X_t \in \cdot \;}{\tau_b > t}) = (-\infty , \min (x, b(t))]
\end{align*}
Since $(X_t)_{t\geq 0}$ has non-increasing paths and $t \in \supp (\PR{\tau_b \in \cdot\;})$ we have
\begin{align*}
\mu ( \{x \geq b(t) : \PRx x {\tau_b > 0}>0 \})>0.
\end{align*}
This implies that
\begin{align*}
\sup \supp (\PRbx \mu {X_t \in \cdot \;}{\tau_b > t}) = b(t).
\end{align*}
Therefore, under the assumptions that \ref{existence: diffusitivy assumption} and ((a) or (b)) are fulfilled, and hence, by Theorem~\ref{uniqueness}, the boundary function $b$ with $\tau_b \overset{\text d} = \xi$ is unique on $I^\xi$.
\end{proof}

\section{Conditions for diffusions on an interval}\label{conditions for diffusions}

In this section we establish conditions under which a diffusion process in an interval, which satisfies a stochastic differential equation up to an explosion time, fulfills the assumptions required for existence and uniqueness of solutions in the inverse first-passage time problem. The proof of Theorem~\ref{diffusions result} is to be found at the end of the section. At first, we will collect the essential steps in preliminary statements. Let $(X_t)_{t\geq 0}$ be a diffusion on an interval $E$ according to Definition~\ref{definition diffusion}.

\begin{proposition}\label{first-passage time representation}
Assume that $R\notin E$. Further, assume that $\sigma \in C^1 ((L,R))$, $\sigma >0$ and that $\beta$ is locally bounded on $(L,R)$. Let $x\in E$. Let $b:[0,\infty) \to [-\infty , \infty ]$ be a boundary function. It holds that
\begin{align*}
\PRx x{\tau_b = \tau_b'}= 1.
\end{align*}
\end{proposition}

In the case of Brownian motion the following statement was proved in Proposition~6.1 in \cite{Eks16}.

\begin{proposition}\label{first-passage time distribution lemma}
Let $b:[0,\infty) \to \overline E$ be a boundary function. Assume that $\sigma \in C^1 ((L,R))$, $\sigma >0$ and that $\beta$ is locally bounded on $(L,R)$. Let $\mu$ be a probability measure on $E$. Assume that $\PRx{\mu}{\tau_b >0} >0$ and that $\PRx{\mu}{X_t\in \cdot}$ is diffuse for every $t>0$. Then
\begin{align*}
\supp ( \PRbx{x}{X_t \in \cdot \;}{\tau_b >t}) = [L,b(t)]
\end{align*}
for every $t< \inf\{s > 0 : b(s) = L\}$.
\end{proposition}

The first step towards Proposition~\ref{first-passage time representation} will be the following.

\begin{lemma}\label{first-passage time representation lemma}
Let $b:[0,\infty) \to \overline E$ be a boundary function. Assume that $\sigma \in C^1 ((L,R))$, $\sigma >0$ and that $\beta$ is locally bounded. Then for $x\in (L,R)$ we have
\begin{align*}
\PRx{x}{\tau_b  < \tau_b' \wedge S} = 0,
\end{align*}
where $S \coloneqq \lim_{n\to\infty} S_n$, where $S_n$ are defined in Definition~\ref{definition diffusion} (iii).
\end{lemma}

The idea is to reduce the situation to Brownian motion and use the fact that for Brownian motion the desired statements are already known. For example, Proposition~2 in \cite{chen2011} and Lemma~6.2 in \cite{Eks16} prove that for Brownian motion it holds $\tau_b = \tau_b'$ almost surely. For this we follow the idea of \cite{chen2011} from Proposition 2 therein. We first scale the process in the spatial coordinate as in (4.2) in \cite{chen2011} and then change the measure by using the Girsanov theorem. 

For simplicity we assume that $\sigma \in C^1 ((L,R))$, $\sigma >0$ and that $\beta$ is locally bounded. Let $c\in (L,R)$ be fixed and for $x\in (L,R)$ define
\begin{align}\label{chen scaling}
f(x) \coloneqq \int_c^x \frac{1}{\sigma (z)}\intdiff z.
\end{align}
We have $f\in C^2(L,R)$ and that $f$ is strictly increasing and invertible. Let $n\in \N$. Under $\Pr_x$ the process $(X_{t\wedge S_n})_{t\geq 0}$ is a semimartingale and due to the It\^o formula it follows that
\begin{align*}
f(X_{t\wedge S_n}) = f(X_0) + \int_0^{t\wedge S_n} \left( \frac{\beta (X_s)}{\sigma(X_s)}  - \frac{1}{2} \sigma' (X_s) \right) \intdiff s + B_{t\wedge S_n}.
\end{align*}
This means that the process given by $Y_t \coloneqq f(X_t)$ fulfills
\begin{align*}
Y_{t\wedge S_n} = Y_0 + \int_0^{t\wedge S_n} \tilde{\beta}(Y_s) \intdiff s + B_{t\wedge S_n},
\end{align*}
where
\begin{align*}
\tilde \beta (y) \coloneqq \frac{\beta (f^{-1} (y))}{\sigma (f^{-1} (y))} - \frac{1}{2} \sigma' (f^{-1}(y)).
\end{align*}
Note that $\tilde{\beta}(Y_t)$ is uniformly bounded in $t \leq S_n$. Let $T>0$ be fixed. Then
\begin{align*}
U_t \coloneqq \exp\left( - \int_0^{ t\wedge T\wedge  S_n} \tilde{\beta}(Y_s) \intdiff B_s - \frac{1}{2} \int_0^{ t\wedge T\wedge S_n } (\tilde \beta (Y_s))^2 \intdiff s \right)
\end{align*}
defines a uniformly integrable positive martingale. By Girsanov the measure
\begin{align}\label{chen scaling girsanov}
\intdiff \tilde{\Pr}_x^{n,T} \coloneqq U_T \intdiff \Pr_x 
\end{align}
defined on $\mathcal{F}_{T\wedge S_n}$ is equivalent to $\Pr_x$ on $\mathcal{F}_{T\wedge S_n}$ and $(Y_{ t\wedge T\wedge S_n})_{t\geq 0}$ is a local martingale. Since its quadratic variation is $(t \wedge T\wedge S_n)_{t\geq 0}$, Lévy's characterization of Brownian motion shows that $(Y_{ t\wedge T\wedge S_n})_{t\geq 0}$ is a Brownian motion stopped at $  T\wedge S_n$.

\begin{proof}[Proof of Lemma~\ref{first-passage time representation lemma}]
Recall $f$ from \eqref{chen scaling} and consider $Y_t = f(X_t)$. Define $\tilde{b} \coloneqq f(b)$, where we allow $f(R) \in (-\infty , \infty ]$ and $f(L) \in [-\infty , \infty )$. Note that
\begin{align*}
\tau_{b} = \inf\{t> 0 : Y_t \geq \tilde{b}(t)\},\quad \tau_b ' = \inf\{t> 0 : Y_t > \tilde{b}(t)\}.
\end{align*}
Since under $\tilde{\Pr}_x^{n,T}$ from \eqref{chen scaling girsanov} the stopped process $(Y_{ t\wedge T\wedge S_n})_{t\ge 0 }$ is a stopped Brownian motion we have that
\begin{align*}
\tilde{\Pr}_x^{n,T}\left( \tau_b < \tau_b' \wedge T \wedge S_n \right) = 0.
\end{align*}
Due to the equivalence of the measures $\tilde{\Pr}_x^{n,t}$ and $\Pr_x$ it follows that
\begin{align*}
\PRx{x}{\tau_b < \tau_b' \wedge T \wedge S_n } = 0.
\end{align*}
By first letting $T\to\infty$ and then $n\to\infty$ it follows that
\begin{align*}
\PRx{x}{\tau_b < \tau_b' \wedge S } = 0.
\end{align*}
This finishes the proof.
\end{proof}

\begin{proof}[Proof of Proposition~\ref{first-passage time representation}]
The idea of this proof is to split the path into suitable excursions away from the lower boundary and to apply Lemma~\ref{first-passage time representation lemma} for every excursion. We begin with assumptions by which we do not lose generality in order to reduce the complexity of the boundary involved.

Since $R\notin E$ we can assume that $b$ takes values in $\overline E$. Due to the a.s. convergence of $\lim_{s\searrow 0} \tau_{b\vert_s} = \tau_b$ and $\lim_{s\searrow 0} \tau_{b\vert_s} ' = \tau_b '$ we can assume that there is $s>0$ such that $b(t) = R$ for $t < s$. Furthermore, we have $\tau_b \leq \tau_b' \leq \inf\{t> 0 : b(t) = L\}$. If $\tau_b < \tau_b' $ we have consequently $\tau_b < \inf\{t> 0 : b(t) = L\}$. We will therefore assume that there is $u < \inf\{t> 0 : b(t) = L\}$ such that $b(t) = R$ for all $t> u$ and $b(t) > L$ for all $t\in [s,u]$. Thus, using the lower semicontinuity, we now treat the case that $b(t) = R$ for $t\notin [s,u]$ and
\begin{align*}
L < \inf_{t\in [s,u]} b(t) \leq \sup_{t\in [s,u]} b(t) \leq R.
\end{align*}
For $x\in (L,R)$ let us define
\begin{align*}
T_x \coloneqq \inf\{t\geq 0 : X_t \leq x\}.
\end{align*}
Let $x_\ell, x_r \in \{\ell_n : n\in \N\}$ such that $L< x_\ell < x_r < \inf_{t\in [s,u]} b(t)$. For $k\in \N$ let us inductively define $\rho_0 \coloneqq 0$,
\begin{align*}
\lambda_{k} &\coloneqq \inf\{t\geq \rho_{k-1}: X_t \leq x_\ell \}
\intertext{and}
\rho_k &\coloneqq \inf\{t\geq \sigma_k : X_t \geq x_r \}.
\end{align*}
Since $X_t \leq x_r < \inf_{t\in [s,u]} b(t)$ for all $t\in [\lambda_k , \rho_k]$ and all $k\in \N$ we have that, if $\tau_b \neq \tau_b'$, then there is $k\in\N_0$ such that
\begin{align*}
\rho_k < \infty \quad \text{and}\quad \tau_b \in (\rho_k , \lambda_k ) \quad\text{and}\quad \tau_b < \tau_b'.
\end{align*}
For $x\in E$ we have that
\begin{align*}
&\PRx{x}{\tau_b \in (\rho_k , \lambda_k ), \tau_b < \tau_b' ,\rho_k < \infty}\\
&= \EWx{x}{\PRx{X_{\rho_k}}{ \tau_{b^\theta} \in (0, T_{x_\ell}) , \tau_{b^\theta} < \tau_{b^\theta} ' }_{\theta={\rho_k}}\IF{\rho_k < \infty }\IF{\rho_k < \tau_b }},
\end{align*}
where $b^\theta (t) = b(\theta +t )$. For the moment fix $\theta \geq 0$  and note that $b^\theta$ is a boundary function taking values in $\overline E$. Since $(X_t)_{t\geq 0}$ has continuous paths and $R\notin E$ we have $\lim_{n\to\infty}\tau_{r_n} = \infty$, where $(r_n)_{n\in\N}$ is the sequence from Definition~\ref{definition diffusion} (iii). Therefore we have that $T_{x_\ell} \leq S = \lim_{n\to\infty} S_n $ and therefore for $z \in (L,R)$, by using Lemma~\ref{first-passage time representation lemma}, it holds
\begin{align*}
\PRx{z}{ \tau_{b^\theta} \in (0, T_{x_\ell}) , \tau_{b^\theta} < \tau_{b^\theta} ' } \leq \PRx{z}{ \tau_{b^\theta } < \tau_{b^\theta } ' \wedge S } = 0.
\end{align*}
Using this for $z = X_{\rho_k} \in \{X_0, x_r\} \subset (L,R)$ and plugging it back into the expectation above we have therefore
\begin{align*}
\PRx{x}{\tau_b \in (\rho_k , \lambda_k ), \tau_b < \tau_b' ,\rho_k < \infty} = 0.
\end{align*}
This implies
\begin{align*}
\PRx{x}{\tau_b \neq \tau_b'} \leq \sum_{k\in\N_0} \PRx{x}{\tau_b \in (\rho_k , \lambda_k ), \tau_b < \tau_b' ,\rho_k < \infty} =0.
\end{align*}
This finishes the proof of the first part.
\end{proof}

\begin{proof}[Proof of Proposition~\ref{first-passage time distribution lemma}]
Let $t < \inf\{s>0 : b(s) =L\}$ and $z\in (L , b(t))$ and $\varepsilon >0$ such that $[z-\varepsilon , z+\varepsilon ] \subseteq  (L , b(t))$. Since $\PRx{\mu}{\tau_b > 0}>0$ there is $s\in (0,t)$ with $\PRx{\mu}{\tau_b > s}>0$. And since $\PRx{\mu}{X_s\in\cdot }$ is diffuse there is $y \in (L,b(s))$ such that $y \in \supp( \PRx{\mu}{X_s \in \cdot\; , \tau_b > s})$. Set $b^s (t) = b(s +t )$. Due to
\begin{align*}
&\PRx{\mu}{ X_t \in ( z-\varepsilon , z+\varepsilon ) , \tau_b > t  }\\
&\geq \int_{\R} \PRx{u}{X_{t-s} \in ( z-\varepsilon , z+\varepsilon ) , \tau_{b^s} > t-s} \PRx{\mu}{X_s \in \intdiff u , \tau_b > s}
\end{align*}
it suffices to show that we have
\begin{align}\label{proof:first-passage time distribution lemma 1}
\PRx{u}{X_{t-s} \in ( z-\varepsilon , z+\varepsilon ) , \tau_{b^s} > t-s}>0
\end{align}
for $u \in U$, where $U \subseteq (L,b(s))$ is a neighborhood of $y$. For this, let $n\in\N$ be large enough such that
\begin{align*}
y\in (\ell_n ,r_n), \quad (z-\varepsilon , z+\varepsilon ) \subseteq (\ell_n , r_n).
\end{align*}
Further, choose $T > t-s$. Recall $f$ from \eqref{chen scaling}. Let $a:[0,t-s]\to\R$ be a continuous function and $\delta >0$ such that
\begin{itemize}
\item $a(0) = f(y)$ and $(a(r) - \delta , a(r) + \delta ) \subseteq (f(\ell_n ) , f(r_n))$ for all $r\in [0,t-s]$,
\item $(a(t-s) - \delta , a (t-s) + \delta ) \subseteq (f(z-\varepsilon) , f(z+\varepsilon ) )$,
\item $f(a(r) + \delta ) < f(b(s+r))$ for all $r\in [0,t-s]$.
\end{itemize}
An explicit construction of the function $a$ can be made as in Lemma~2.3.6 of \cite{k2022}. The last point is possible since $b(sr) >L$ for all $r\in [0,t-s]$ and $b$ is lower semicontinuous. Recall $Y_r = f(X_r)$ and $\tilde{\Pr}_{u}^{n,T}$ from \eqref{chen scaling girsanov} and that under $\tilde{\Pr}_{u}^{n,T}$ the stopped process $(Y_{ r\wedge T\wedge S_n})_{r\ge 0 }$ is a stopped Brownian motion. Note that it happens with positive probability that a Brownian motion started at $f(u)\in (f(y)-\delta, f(y)+\delta)$ stays up to time $t-s$ in a tube which follows a continuous function, for instance see Theorem 38 of \cite{Freedman1983}. This argument was already used in Proposition 3.1 of \cite{Eks16} in case of Brownian motion.
This leads to
\begin{align*}
&\tilde{\Pr}_{u}^{n,T}\left( Y_{t-s} \in ( f(z-\varepsilon) , f(z+\varepsilon) ) , \tau_{b^s} > t-s,  S_n > t-s  \right) \\
&\geq \tilde{\Pr}_{u}^{n,T}\left( \abs{ Y_r -  a(r)} < \delta \;\forall r\in [0,t-s]\right)> 0.
\end{align*}
Due to the equivalence of the measures $\tilde{\Pr}_x^{n,t}$ and $\Pr_x$ we obtain that \eqref{proof:first-passage time distribution lemma 1} is true. All in all, since $z$ was arbitrary, this means that 
\begin{align*}
\supp (\PRbx{x}{X_t \in\cdot\;}{\tau_b >t}) = [L,b(t)].
\end{align*}
This finishes the proof of the statement.
\end{proof}

\begin{proof}[Proof of Theorem~\ref{diffusions result}]
By our definition of a diffusion on an interval we already assumed that the process has continuous paths and is a strong Markov process. This gives \ref{existence: continuity assumption} and \ref{uniqueness: markov assumption}. Since $(X_t)_{t\geq 0}$ is a strong Markov process and has continuous paths, the transition probabilities preserve the usual stochastic order, since paths that started from different positions can be let run together after they have met, for details see Lemma~A.6.1 in \cite{k2022}. Hence we have \ref{uniqueness: order preservation assumption}. Moreover, due to the assumptions on the coefficients and that $R\notin E$ we can apply Proposition~\ref{first-passage time representation} and obtain \ref{existence: first-passage time assumption}. 

Now assume that \ref{existence: diffusitivy assumption} holds. Let $t\in (0,t^\xi)$ and assume that $\tau_b \overset{\text d} = \xi$. Since $\PRx{\mu}{\tau_b >t}>0$ we have that $t< \inf\{s>0 : b(s) =L\}$. Therefore, Proposotion~\ref{first-passage time distribution lemma} implies that
\begin{align*}
\sup \supp (\PRx{\mu}{X_t\in\cdot\;,\tau_b> t}) = b(t).
\end{align*}
This gives \ref{uniqueness: support assumption}, and therefore there exists a boundary function $b$ with $\tau_b \overset{\text d} = \xi$ is unique on $(0,t^\xi)$.
\end{proof}

\appendix

\section{Certain discontinuity sets for arbitrary functions}

The following statement follows from Theorem~6 in \cite{young1907}, but for completeness we give an own proof, which makes use of probabilistic arguments.

\begin{lemma}\label{discontinuity lemma}
Let $b:[0,\infty] \to [-\infty,\infty]$ be an arbitrary function. Then the set
\begin{align*}
\left\{ t\in (0,\infty) : \max\left(\liminf_{s\nearrow t} b(s) , \liminf_{s\searrow t} b(s)\right) > b(t) \right\}
\end{align*}
is countable.
\end{lemma}
\begin{proof}
We only consider the set
\begin{align*}
S_b \coloneqq \left\{ t\in (0,\infty) : \liminf_{s\nearrow t} b(s) > b(t) \right\}
\end{align*}
since then the statement follows for the remaining points by consideration of the map $(0,\infty)\ni t\mapsto b ( 1/ t)$.

Further, let
\begin{align*}
\varphi : [-\infty, \infty] \to [-1,1] ,\; \varphi (x) \coloneqq   \frac{x}{1+ \abs{x}} \If_{\R}(x) + \sgn (x) \If_{\{-\infty , \infty\}}(x)
\end{align*}
and set $\tilde{b}(t) \coloneqq \varphi (b(t))$. Since $S_b \subseteq S_{\tilde{b}}$ we can assume that $b$ takes values in $[-1,1]$.

The function defined by
\begin{align*}
b^* (t) \coloneqq \min \left( \liminf_{s\to t} b(s) , b(t)\right)
\end{align*}
is lower semicontinuous and it holds $S_b \subseteq S_{b^*}$. Thus without loss of generality we can assume that $b$ is a lower semicontinuous function.

Let $t\in S_b$. Then there is $\varepsilon >0$ such that there exists $\delta>0$ with
\begin{align*}
b(s) \geq b(t) + \varepsilon \quad \forall s \in (t-\delta , t).
\end{align*}
Let $(B_t)_{t\geq 0}$ be a Brownian motion starting from a deterministic point $B_0 = x < -1$. For a function $f:[0,\infty]\to \R$ define
\begin{align*}
\tau_f \coloneqq \inf\{s>0 : B_s \geq f(s)\}.
\end{align*}
Let $K\coloneqq b(t) +\varepsilon$. Then we have
\begin{align*}
\PR{\tau_b = t} &\geq \PR{\tau_{-1}> t-\delta, \tau_{K}> t, B_t \in (b(t),b(t)+\varepsilon)}\\
&= \PR{ \sup_{s\in [0,t-\delta]} B_s < -1 , \sup_{s\in[t-\delta ,t]} B_s < K, B_t \in (b(t),b(t)+\varepsilon) } >0.
\end{align*}
Therefore, we have
\begin{align*}
S_b \subseteq \{t\in (0,\infty) : \PR{\tau_b = t} > 0 \},
\end{align*}
where the right-hand side is a countable set. This finishes the proof.
\end{proof}
%
%

\section*{Acknowledgments}

Alexander Klump gratefully acknowledges the support of a postdoctoral fellowship from the German Academic Exchange Service (DAAD) and the hospitality of the Institute of Mathematics and Informatics of the Bulgarian Academy of Sciences.

Mladen Savov acknowledges - "This study is financed by the European Union-NextGenerationEU, through the National Recovery and Resilience Plan of the Republic of Bulgaria, project No BG-RRP-2.004-0008"

\printbibliography

\end{document}